\documentclass[a4paper,10pt]{amsart}
\usepackage{amssymb}
\usepackage{amscd}

\numberwithin{equation}{section}

\let\ep=\varepsilon
\def\n{|\!|}
\def\cinf{{\mathcal{C}^\infty}}

\def\cb{\mathbb{C}}
\def\nb{\mathbb{N}}

\def\pb{\mathbb{P}}
\def\rb{\mathbb{R}}
\def\zbb{\mathbb{Z}}
\def\oc{{\mathcal{O}}}
\def\ot{\otimes}
\def\det{\text{det}}
\def\xb{\overline{x}}
\def\zb{\overline{z}}
\def\phi{\varphi}

\newtheorem{theo}{Theorem}
\newtheorem{prop}{Proposition}

\newtheorem{lemma}{Lemma}

\begin{document}
\title{Computations of Bott-Chern classes on $\mathbb{P}(E)$}
\author{Christophe Mourougane}
\address{Institut de Math{\'e}matiques de Jussieu / Plateau 7D /
175, rue du  Chevaleret /
 75013 Paris}
\email{christophe.mourougane@math.jussieu.fr}
\maketitle

{\textbf{Abstract.}
\textit{ We compute the Bott-Chern classes of the metric Euler
  sequence describing the relative tangent bundle of the variety  
$\mathbb{P} (E)$
  of hyperplans of a holomorphic hermitian vector bundle $(E,h)$
 on a complex manifold.
We give applications to the construction of the arithmetic characteristic
classes of an arithmetic vector bundle $\overline{E}$ 
and to the computation of the height of $\mathbb{P}(\overline{E})$
 with respect to the tautological quotient bundle 
$\mathcal{O}_{\overline{E}} (1)$.}
~\footnote{Key words : Bott-Chern secondary classes, metric relative 
Euler sequence, arithmetic characteristic classes.\\ MSC : 32L10, 14G40}

\section*{Introduction}
\label{sec:intro}
In the whole work, $X$ will be a smooth complex analytic manifold
of dimension~$n$. For a vector bundle $F\to X$ on~$X$ we denote by
$A^{p,q}(X,F)$ the space of smooth $F$-valued differential forms on
$X$ of type $(p,q)$. 
On a holomorphic vector bundle $E\to X$
endowed with a hermitian metric $h$, there exists a unique connection
$\nabla =\nabla _{E,h}$ compatible with both the holomorphic and the
hermitian structures. It is called the Chern connection of $(E,h)$.
 Its curvature $\frac{i}{2\pi}\nabla ^2$ is multiplication 
by an endomorphism-valued $2$-form that we will denote by
$\Theta(E,h)\in A^{1,1}(X,Herm(E))$. 
The associated Chern forms are defined by 
$$\det (I+t A)=\sum_{d=0}^{+\infty} t^d \det_d (A) \text{ and }
  c_d(E,h):= \det_d(\Theta(E,h))\in A_{d,d}(X,\cb). $$
The form $c_d(E,h)$ is closed and represents the Chern class $c_d(E)$ 
of $E$ in $H^{2d}(X,\rb )$ through De Rham isomorphism.
The Chern class polynomial $c_t(E):=\sum_{d=0}^{+\infty} t^d c_d(E)$ is
multiplicative on exact sequences. Bott-Chern secondary classes are classes in 
$$\widetilde{A^{d,d}}(X):=\frac{A^{d,d}(X,\cb)}{Im d'+Im d''}$$
 which functorially represent the default of multiplicativity of the
 Chern form polynomial $c_t(E,h):=\sum_{d=0}^{+\infty} t^d c_d(E,h)$
 on short exact sequences of hermitian vector bundles
 (see~\cite{bis-gs} theorem 1.29 and~\cite{gs} theorem 1.2.2)~:
 \begin{itemize}
 \item For any short exact sequence $(\mathcal S ) =(0\to  S\to E\to Q\to
   0)$, and any choice of metrics $h=(h_E,h_S,h_Q)$, 
$\widetilde{c_{d+1}}(\mathcal S,h) \in\widetilde{A^{d,d}}(X)$,
$$c_t(E,h_E)-c_t(S,h_S)c_t(Q,h_Q)=
-\frac{it}{2\pi}d'd''\widetilde{c_t}(\mathcal S,h)$$
 where $\widetilde{c_t}(\mathcal S,h)=\sum_{d=0}^{+\infty} t^d
   \widetilde{c_{d+1}}(\mathcal S,h)$. 
\item If $(\mathcal S,h)$ is metrically split, 
then $\widetilde{c_t}(\mathcal S,h)=0$.
\item For any holomorphic map $f~:X\to Y$ of complex analytic manifolds and any
  metric short exact sequence $(\mathcal S,h)$ over $Y$, 
$\widetilde{c_t}(f^\star\mathcal S,f^\star h)
=f^\star\widetilde{c_t}(\mathcal S,h)$.
 \end{itemize}

\bigskip
They were introduced by Bott and Chern in their study of the 
distribution of the values of holomorphic sections of hermitian vector
bundles~\cite{bc}. They were used by Donaldson for  defining
 a functional on the space of hermitian metrics on a given holomorphic
 vector bundle in order to find Hermite-Einstein metrics~\cite{do}. They were
 given an axiomatic definition and served
 as one kind of secondary objects (together with Green currents and
 analytic torsion forms) in the study by Bismut Gillet and
 Soul{\'e} of holomorphic determinant bundle~\cite{bis-gs}. They enter the very
 definition of arithmetic characteristic classes set by Gillet and
 Soul{\'e}~\cite{gs}.

\bigskip
Bott-Chern classes were computed in few cases, mainly when $(E,h)$ is
 assumed to be flat. See the work of Gillet and Soul{\'e} in the the case of
 projective spaces~\cite{gs} of Maillot in the case of
 Grassmanians ~\cite{mai} and of Tamvakis in the case of other flag
 manifolds~\cite{tam}. In most of computed cases, 
 Bott-Chern classes are made of closed forms.

  \bigskip 
We deal with a relative situation.
 Let $E\to X$ be a holomorphic vector bundle of rank~$r$ over
$X$ . Consider  $\pi :\pb (E)\to X$, the variety of hyperplans of
$E$. The differential of the~quotient map 
$E^\star-X\times \{0\}\to \pb(E)$ over $X$ gives rise to the relative
Euler sequence~$(\Sigma )$ on  $\pb (E)$~:
\begin{eqnarray*}
  0\to \oc _E(-1)\to\pi^\star E^\star\to T\to 0 
\end{eqnarray*}
where $T$ denotes the twisted relative tangent bundle 
$T_{\pb (E)/X}\ot \oc_E(-1)$.
The choice of a~hermitian metric $h$ on $E$ enables to 
endow all the bundles in the
sequence with a natural hermitian metric.
Our aim is to compute Bott-Chern forms (i.e. representatives of 
Bott-Chern secondary classes) for this metric relative Euler sequence
$(\Sigma,h)$ in terms of the curvature form of $(E,h)$.
The case where $X$ is a point is dealt with by Gillet and Soul{\'e}
in~\cite{gs}. 

\bigskip
We denote by $\alpha$ the curvature form $\Theta(\oc_E(1),h)$ of
 $(\oc_E(1),h)$. For $q\in\nb$,
we consider the forms $\Theta ^q a^\star$ 
on $\pb (E)$ given  at the point $(x,[a^\star ])$ of $\pb (E)$, 
($a^\star\in E^\star_x$) by
$$\displaystyle\Theta ^q a^\star=\frac{\langle\pi^\star\Theta(E
  ^\star,h)^q a ^\star,a^\star\rangle}{\n a ^\star\n^2}.$$ 
The product is taken with respect to the wedge product in the form
part and the composition in the endomorphism part.
The generating function $\Theta_t$ for those
forms is $\sum_{d=0}^{+\infty}(-t)^d (\Theta ^d a ^\star)$.
We will also need the Harmonic numbers
$\mathcal{H}_p=\sum_1^p\frac{1}{i}$  and their generating function
$H(X):=\sum_{p=1}^{+\infty}\mathcal{H}_p X^{p}$.
We prove

\bigskip
\textbf{Main theorem.}
\textit{ The Bott-Chern secondary class polynomial
$\widetilde{c_t}(\Sigma,h)$ for $(\Sigma,h)$ is represented by the form
$$-H\left(t\alpha\Theta_t
+1-\Theta_t
+t\Theta_t\frac{i}{2\pi}d'd''\ln\Theta_t \right)
\Theta_t\pi ^\star c_t(E ^\star,h).$$
}

This form is in general not closed.
The first part of the proof consists in finding an expression in
coordinates for $\widetilde{c_{d+1}}(\Sigma,h)$.  Computations are
simplified by the choice of a normal frame at a given point
$(x_0,[a_0^\star ])$ of $\pb (E)$. The combinatorial analysis is
rather intricate.
The second part is devoted to identify the previously found
expressions in terms of globally defined forms. Finally, the
computation of the Bott-Chern forms generating polynomial
enables to express the result in a concise way, for Bott-Chern forms
follow a recursion formula.

\bigskip Our second aim is to give some applications in the setting of
Arakelov geometry. 
Let $\mathcal{X}$ be an arithmetic variety built on a scheme $\chi$
defined on  the ring of integers of some number field. Let $\overline{E}$ be
 an arithmetic vector bundle on $\mathcal{X}$. We will 
write $X$ for the smooth variety $\chi_\cb(\cb)$ and $(E,h)$ 
for the induced hermitian vector bundle on $X$.
We give in section~\ref{arithmetic}
an application of our computations to the construction of
arithmetic Chern classes from arithmetic Segre classes. 
At algebraic level, the total Chern class of $E^\star$ is the inverse of the
total Segre class of $E$. This is also true at the level of forms
(see section~\ref{form}). But
at the arithmetic level, some secondary purely complex analytic
objects have to enter the definition
of Segre classes in order to keep this inverse relation true in the
arithmetic Chow group of the base scheme. This is due to the fact that
$\overline{E}^\star \ot \oc_{\overline{E}} (1)$
 having a nowhere vanishing section
has nevertheless a non-zero top Chern form. 
Evaluating some fiber integrals of previously studied Bott-Chern
forms, we give explicit expressions for the secondary classes involved.
This answers a question raised by Bost and Soul{\'e}~\cite{bost}.
This in turn enables to compute the height of $\pb(\overline{E})$
 with respect to $\oc_{\overline{E}} (1)$ in term of the top arithmetic
Segre class and a secondary term from complex cohomology. The
qualitative output is that complex characteristic classes of $E$
are the only complex datas of $(E,h)$ 
that accounts for the measure of the complexity of
 $\pb(\overline{E})$ given by its height. 

No wonder that such computations at least in the full relative case
(see proposition~\ref{prop:full}) can be done for other characteristic
classes than the Chern classes on other flag varieties
of~$\overline{E}$ to get expressions for their height.

\thanks{I thank Vincent Maillot for his interest in my work and for his
suggestion to study the arithmetic applications.}

\tableofcontents

\section{Bott-Chern forms}
\label{sec:bott}
General formulas for computing $\tilde c(\mathcal S,h)$ with induced
and quotient metrics were given by Bott
and Chern ~\cite{bc}. For sake of completeness, we outline their
method.

The commutator of $End(E)$-valued forms $\alpha$ and $\beta$ is
defined to be
$[\alpha,\beta]:= \alpha\beta-(-1)^{\deg\alpha
  \deg\beta}\beta\alpha.$ 

Consider the sequence 
$0\to S\stackrel{\iota}\to E\stackrel{p}\to Q\to O$ 
endowed with metric
constructed from a hermitian metric $h$ on $E$. Denote by $\nabla$
(resp. $\nabla_S$, $\nabla _Q$) the
Chern connection of $(E,h)$ (resp. $(S,h_{|S})$, $(Q,h_{|Q})$).
Consider the family of connections on $E$
$$\nabla_u:= \nabla +(u-1)P_Q\nabla P_S$$
where $P_S=\iota\iota ^\star$ (resp. $P_Q=p^\star p$) denotes the orthogonal
projection of $E$ onto $\iota(S)$ (resp. $\iota(S)^\perp$). 
The choice of a local holomorphic frame for $E$  enables to 
express locally the connection $\nabla_u$ as $d+A_u$ for some matrix
valued $(1,0)$-form $A_u$. The curvature $\Theta(\nabla_u)$ of
$\nabla_u$ is then given by the matrix valued $(1,1)$-form 
$\frac{i}{2\pi}(dA_u+A_u\wedge A_u)$. Recall that the connection
$\nabla _u$ extends to $End(E)$-valued forms by $\nabla _u
(\alpha):=[\nabla _u ,\alpha]$. The
formula for $\tilde c(\Sigma ,h)$ relies on the following identities
\begin{eqnarray*}
\nabla_u (\Theta(\nabla_u))&=&[\nabla _u , \Theta(\nabla_u)]= [\nabla _u
  ,\nabla _u ^2]=0  \text{        Bianchi formula}\\
&=& d(\Theta(\nabla_u)) + [A_u, \Theta(\nabla_u)]\\
\nabla_u (P_S)&=&[\nabla_u,P_S]=u\frac{d}{du}\nabla_u\text{ got from }
  P_S\nabla' P_Q=P_Q\nabla''P_S=0\\ 
&=& d'P_S+[A_u,P_S].\\
\end{eqnarray*}
That is
\begin{eqnarray*}
d'\Theta(\nabla_u)=d\Theta(\nabla_u)&=&-[A_u, \Theta(\nabla_u)]\\
d'P_S&=& -[A_u,P_S]+u\frac{d}{du}\nabla_u.
\end{eqnarray*}

Consider the polarization $Det_{k}$ of $\det_k$, that is the
symmetric $k$-linear form on $M_r(\cb)$ whose restriction on  the
small diagonal is $\det_k$. Note by $Det_k(A;B)=k Det_k(A,A,\cdots,A,B)$.
The differential version of the $Gl(r,\cb)$-invariance of $\det_k$
shows that the contribution of the commutators $-[A_u, \cdot ]$
vanishes. This leads to 
 \begin{eqnarray*}
d'\left( Det_k
  (\Theta(\nabla_u);P_S)\right)=u Det_k\left(\Theta(\nabla_u);
\frac{d}{du}\nabla_u \right).
\end{eqnarray*}
Now, 
\begin{eqnarray*}
\frac{i}{2\pi}d\left(\frac{d}{du}\nabla_u \right)&=&\frac{i}{2\pi}
\frac{d}{du} d(\nabla_u)
=\frac{i}{2\pi}\frac{d}{du} d(A_u) 
=\frac{d}{du} \left( \Theta(\nabla_u)-\frac{i}{2\pi} A_u\wedge
  A_u \right)\\
&=&-\frac{i}{2\pi}[A_u ,\frac{d}{du}A_u] +\frac{d}{du} \left(
  \Theta(\nabla_u)\right).
\end{eqnarray*}
The $Gl(r,\cb)$-invariance of $\det_i$ leads to
\begin{eqnarray}
\nonumber
-\frac{i}{2\pi}d'd'' Det_k (\Theta(\nabla_u);P_S)&=&
\frac{i}{2\pi}d\left(u Det_k\left(\Theta(\nabla_u);
\frac{d}{du}\nabla_u \right)\right)\\
&=& u\frac{d}{du}\det_k(\Theta(\nabla_u)). \label{diff}
\end{eqnarray}

One then checks in a frame adapted to the $\cinf$ splitting $E\simeq
S\oplus Q$ given by $\iota^\star\oplus p$, that 
$$\Theta(\nabla_u)=\left| 
  \begin{array}{cc}
\mbox{}&\mbox{}\\
(1-u)\Theta_S + u\iota^\star\Theta_E \iota & \iota^\star\Theta_E p^\star\\
\mbox{}&\mbox{}\\
up\Theta_E \iota& (1-u)\Theta_Q + up\Theta_E p^\star\\
\mbox{}&\mbox{}
  \end{array}
\right|$$
Hence integrating equation~{(\ref{diff})} between $0$ and $1$, we get 
\begin{eqnarray}\label{1}
c_k(E,h)-c_k(S\oplus Q,\nabla_S\oplus\nabla_Q)=-\frac{i}{2\pi} d'd''
\int_0^1\frac{\Phi_k(u)-\Phi_k(0)}{u}du
\end{eqnarray}
 where 
\begin{eqnarray*}
  \Phi _k(u)&=&Det_k(\Theta(\nabla_u);P_S)=\text{coeff}_\lambda
det_k(\Theta(\nabla_u)+\lambda P_S).
\end{eqnarray*}
If moreover, the sub-bundle $S$ is of rank $1$, then 
\begin{eqnarray}\label{2}
 \Phi _{d+1}(u)=\det _{d}\left(\left(1-u\right)\Theta_Q + up\Theta_E
  p^\star  \right). \end{eqnarray}

\bigskip\bigskip
\section{Curvature computations}
\label{sec:curvature}

We will compute the curvature of $T$ and of $\oc _E(1)$ 
at a point $(x_0,[a^\star_0])$ of $\pb(E)$ ($a^\star_0~\in~E^\star$) in
a well-chosen frame. 

\bigskip
We now recall the formula for the curvature of a quotient bundle.
Consider the sequence 
$0\to S\stackrel{\iota}\to E\stackrel{p}\to Q\to O$ endowed with metrics
constructed from a hermitian metric on $E$. 
All the bundles we will consider, including the endomorphism bundles
inherit a metric and a Chern connection from those of $E$.
 We will denote with a star the adjoint maps with respect to the given metric. 
Notice first that since $pp^\star=Id_Q$ and $\nabla_E$ is compatible with
the metric, $p(\nabla_{Hom(Q,E)} p^\star)=0$. 
Let $s$ be a local smooth section of $Q$.
We get first
$$\nabla_{E}(p^\star s)=(\nabla_{Hom(Q,E)} p^\star)s+p^\star\nabla_Q s =
\iota\iota^\star(\nabla_{Hom(Q,E)}  p^\star)s + p^\star\nabla_Q s $$
and then 
$$p\nabla ^2_E(p^\star s)=p(\nabla_{Hom(S,E)}\iota )\wedge
\iota^\star(\nabla_{Hom(Q,E)} 
  p^\star)s+\nabla_Q^2 s.$$ 
That is, taking the type into account, (${\nabla''}_{Hom(S,E)}\iota=0$)
$$\Theta _Q=p\Theta_E p^\star- \frac{i}{2\pi} p(\nabla'_{Hom(S,E)}\iota)
\wedge\iota ^\star(\nabla''_{Hom(Q,E)} p^\star).$$

\bigskip
We choose local coordinates $x_1,\cdots ,x_n$ around $x_0$, and
a frame $e_1^\star,\cdots ,e_r^\star$ for $E ^\star$ around $x_0$,
normal at $x_0$ 
(i.e. $\langle e_i^\star,e_j^\star\rangle =\delta_{ij}
-2\pi\sum_{\lambda\mu} c_{\lambda\mu ji}x_\lambda\xb_\mu +O(|x|^3)$)
 and such that $a^\star_0=e_1^\star(x_0)$. 
Notice that in this frame at $x_0$, $c_{\lambda\mu ij}
=\overline{c_{\mu\lambda ji}}$, the connection $\nabla_{E^\star}$ is equal
to $d$ and the
curvature $\Theta(E ^\star,h)$ to $i\sum_{\lambda\mu jk}c_{\lambda\mu
  jk}dx_\lambda \wedge d\xb_\mu (e_k^\star)^\star\otimes e_j^\star $. 
Recall the Euler sequence~:
$$0\to\oc_E (-1) \stackrel{\iota}\to \pi^\star E^\star \stackrel{p}\to 
 T\to 0.$$
On an appropriate open set around $(x,a^\star_0)$,
the map $q~:E^\star-X\times \{0\}\to \pb (E)$ is given in coordinates by
$$q(x,\sum a_i
e_i^\star)=(x,[a_1:\cdots:a_r])=(x,\frac{a_2}{a_1},\cdots,\frac{a_r}{a_1}).$$
Hence the map $p$ is given by
$$p\left(x,[\sum_1^r a_i e_i^\star], \sum_1^r b_j e_j^\star\right)=
\left(x,[\sum a_i e_i^\star],\sum_{j\geq 2}
\frac{b_ja_1-b_1a_j}{a_1^2}\frac{\partial}{\partial
  z_j}\right)\ot(\sum a_i e_i^\star).$$ 
where for $2\leq j\leq r$, $z_j:=\frac{a_j}{a_1}$. Here and in the
sequel we will also write $z_1:=\frac{a_1}{a_1}=1$ and $dz_1=0$
 for convenience.
The adjoint of the map $p$ is given by
\begin{eqnarray*}
p^\star\left(\frac{\partial}{\partial z_j}\ot\frac{\sum
  a_i e_i^\star}{a_1}\right)&= &
 e_j^\star-\frac{\langle e_j^\star,\sum a_ie_i^\star\rangle}{\n\sum a_ie_i
    ^\star\n^2} \sum a_ie_i^\star.
\end{eqnarray*}
So, at the point $(x_0,[a ^\star_0])$, the normality of the frame gives
\begin{eqnarray*}
\iota^\star({\nabla''} p^\star)
\frac{\partial}{\partial z_k}\ot\frac{\sum a_i e_i^\star}{a_1} 
&=& \iota^\star{\nabla''}\left( p^\star
\frac{\partial}{\partial z_k}\ot\frac{\sum a_i e_i^\star}{a_1}\right)
-\iota^\star p^\star {\nabla''}
\left(\frac{\partial}{\partial z_k}\ot\frac{\sum a_i e_i^\star}{a_1}\right)\\
 &=& \iota^\star{\nabla''}\left( p^\star
\frac{\partial}{\partial z_k}\ot\frac{\sum a_i e_i^\star}{a_1}\right)
=-d\zb_k\ot\frac{\sum a_i e_i^\star}{a_1}.
\end{eqnarray*}
Now,
\begin{eqnarray*}
  p(\nabla' \iota)\left(\frac{\sum a_i e_i^\star}{a_1}\right)
&=&p\nabla' \left(\iota\frac{\sum a_i e_i^\star}{a_1}\right)
=
\sum_2^r dz_j\frac{\partial}{\partial z_j}\ot\frac{\sum a_i e_i^\star}{a_1}.
\end{eqnarray*}
Summing up the previous computations we infer the formula for the
curvature of $T=T_{\pb(E)/X}\ot\oc_E(-1)$ at $(x_0,[a^\star_0])$
\begin{eqnarray}\label{3}
\Theta(T,h)=\sum _{2\leq j,k\leq r} (c_{jk}+ \frac{i}{2\pi} dz_j\wedge d\zb_k )
  \left(\frac{\partial'}{\partial z_k}\right)^\star\ot
\left(\frac{\partial'}{\partial z_j}\right)
\end{eqnarray}
where 
$\Theta(E ^\star,h)=\sum_{1\leq j,k\leq r}c_{jk}(e_k^\star)^\star \otimes
e_j^\star$, $c_{jk}:=i\sum_{\lambda \mu}
c_{\lambda\mu jk} dx_\lambda \wedge 
d\xb_\mu$ ($c_{jk}=\overline{c_{kj}}$)
 and $\left(\frac{\partial'}{\partial z_j}\right)$ is 
$\left(\frac{\partial}{\partial z_j}\right)\ot\frac{\sum a_i e_i^\star}{a_1}$. 

\bigskip
For later use, we give the formula for the curvature of $(\oc_E(1),h)$.
 Denote by $\Omega$ the positive form defined on the whole fiber of
 $\pi$ over $x_0$ 
 as the Fubini-Study metric of $(\pb (E_{x_0}),h)$~: in coordinates
 $$\Omega:=\frac{i}{2\pi}
\frac{(1+|z|^2)\sum_{j\geq 2} dz_j\wedge d\zb_j-\sum_{i,j\geq 2} \zb_i
  z_j dz_i\wedge d\zb _j}{(1+|z|^2)^2}.$$ 
Then,
\begin{eqnarray*}
\alpha :=\Theta(\oc_E(1),h)= \frac{i}{2\pi} d'd''\log\n e_1 ^\star+
\sum_{j\geq 2}  z_j e_j^\star\n^2 
= \Omega-\frac{\langle\pi^\star\Theta(E ^\star,h)a ^\star,a
 ^\star\rangle}{\n a ^\star\n^2}.
\end{eqnarray*}
This equality is valid on the whole  $\pi ^{-1}(x_0)$.
 At the point $(x_0,[a ^\star_0])$,
it reduces to $\alpha = \frac{i}{2\pi}\sum_{j\geq 2} dz_j\wedge
d\zb_j -c_{11}$.

\bigskip\bigskip
\section{Computations in coordinates}
\label{sec:calcul}
We will compute the Bott-Chern forms of the metric Euler sequence 
at a point $(x_0,[a^\star_0])$ of $\pb(E)$. The results will be given
in terms of locally defined forms.

\bigskip
By the previous results~(\ref{1}),~(\ref{2}) and~(\ref{3}), 
$$\Phi _{d+1}(u)=\det_{d}(\left|c_{jk}+ (1-u)\frac{i}{2\pi} dz_j\wedge
  d\zb_k \right|_{2\leq j,k\leq r})$$
and $$\widetilde{c_{d+1}}=\int_0^1\frac{\Phi_{d+1}(u)-\Phi_{d+1}(0)}{u}du.$$
Our aim is to give an explicit coordinate free expression for the
Bott-Chern forms $\widetilde{c_{d+1}}$.

\bigskip
For $d=0$, $\Phi _1(u)=1$ and $\widetilde{c_1}=0$.

For $d=1$, 
\begin{eqnarray*}
 \Phi _2(u)&=& Trace (\left|c_{jk}+ (1-u)\frac{i}{2\pi} dz_j\wedge
  d\zb_k \right|_{2\leq j,k\leq r}) \\ &=& c_1(\pi ^\star E ^\star) -c_{11}
  +(1-u)\Omega
  \end{eqnarray*}
  and $\widetilde{c_2}=-\Omega$.

 \bigskip 
We need some notations to go along the next computations.
The set of positive integers will be denoted by $\nb_\star$.
For any positive integer $r$, we denote by $\nb_r:=\{1,2,\cdots, r\}$
and $\Sigma_r$ the group of permutations of $\nb_r$. 
The notation $(a_1,a_2,\cdots , a_p)$ will be used for the cycle of
length $p$, $a_1\mapsto a_2\mapsto \cdots \mapsto a_p\mapsto a_1$.

For a finite sequence $B$ in
$\nb_\star^{(\nb_\star)}$ of positive integers, we set  $l(B)=c$ such
that $B$ belongs to $\nb_\star^c$ (\textit{its length})
 and $|B|=\sum_{i=1}^c b_i$
(\textit{its weight}).

A finite sequence $P$ in
$\nb_\star^{(\nb_\star)}$ of positive integers will be called a
\textit{partition} if it is non-decreasing. Its elements will then be called
the parts. If  $[1+\sum _{j=1}^i h_j ,\sum
_{j=1}^{i+1} h_j]$, $(0\leq i\leq q-1, |h(P)|=\sum _{j=1}^q h_j=l(P))$ are the
biggest intervals where the parts are constant
i.e. ~$p_1=p_2=\cdots =p_{h_1}<p_{h_1+1}=\cdots
=p_{h_1+h_2}<p_{h_1+h_2+1} \cdots \cdots p_{h_1+h_2+\cdots h_{q-1}}\\
<~p_{h_1+h_2+\cdots h_{q-1}+1} =~\cdots = p_{h_1+h_2+\cdots h_{q-1}+h_q} 
=p_{l(P)}$ then the sequence $h(P)~=~(h_i)_{1\leq i\leq q}$ will be called
\textit{the height} of the partition $B$. 
 The symbol $h(B)!$ will denote $\prod _{i=1}^q h_i !$.

For a partition $P$ in $\nb_\star^{(\nb_\star)}$ of positive integers,
 a permutation $\sigma$ in $\Sigma_{r}$ will be said to be of
\textit{type $P$} if $r=|P|$ and $\sigma$ can be written as the
 product of $l(P)$ cycles $C_i$ (with disjoint support) 
where $C_i$ is of length $p_i$. 

\bigskip 
We define
\begin{eqnarray*}
c'_i(E ^\star)&:=&\det_{i}
(\left|c_{jk}\right|_{2\leq j,k\leq r}).
\end{eqnarray*}
For a cycle of length $p$, we will need ($p\in\nb_\star,
s\in\nb_\star, s\leq p$, $Q\subset\nb_p$)
$$\Omega_{p,s}:=\frac{1}{p}\sum_{{Q\subset\nb_p\atop\sharp Q=s}}
\Omega_{p,Q} \textrm{ with } \Omega_{p,Q}:=\sum_{2\leq
  i_1,i_2,\cdots,i_p\leq r} \bigwedge_{a=1}^p C_{a,a+1}^{I,Q}
$$
\begin{eqnarray*}\displaystyle
\textrm{ where }&&\left\{
\begin{array}{ccc}
C_{a,a+1}^{I,Q}&=&\left(\frac{i}{2\pi}\right)
dz_{i_a}\wedge d\zb_{i_{a+1}} \textrm{ if } a\in Q\\
&=& c_{i_a i_{a+1}} \textrm{ otherwise }
\end{array}\right.
\end{eqnarray*}
where we have set $i_{p+1}:=i_1$.
 The set 
$Q$ denotes the locations of the $dz\wedge d\zb$ terms. Hence, $s$
is the total number of  $dz\wedge d\zb$ terms occurring in the cycle.
Note that 
$$\Omega_{p,1}=\frac{i}{2\pi}\sum_{2\leq i_1,i_2,\cdots,i_p\leq r}
c_{i_1i_2}c_{i_2i_3}\cdots c_{i_{p-1}i_p}dz_{i_p}\wedge d\zb_{i_1}.$$ 
\bigskip\bigskip

We begin by computing the leading coefficient of our coming formula
for $\Phi$. 
\begin{lemma}
\label{lem:omegap}
 For fixed $d\in\nb_\star$,
$$\sum_{S,P\in\nb_\star^{(\nb_\star)}\atop S\leq P, |P|\leq d, |S|=d} 
\frac{(-1)^{|P|+l(P)}}{l(P)!} \bigwedge_{i=1}^{l(P)}\Omega_{p_i,s_i}
=\Omega ^d.$$
\end{lemma}

\begin{proof}
First notice that 
$\Omega_{p,p}=\frac{1}{p} \Omega_{p,\nb_p} = \frac{1}{p}
(-1)^{p-1}\Omega ^p$. Now,
 \begin{eqnarray*}
\sum_{P,S\in\nb_\star^{(\nb_\star)}\atop S\leq P, |P|\leq d, |S|=d} 
\frac{(-1)^{|P|+l(P)}}{l(P)!}
\bigwedge_{i=1}^{l(P)}\Omega_{p_i,s_i}&=&
\sum_{P,S \in\nb_\star^{(\nb_\star)}\atop S=P,  |P|=d}  
\frac{(-1)^{|P|+l(P)}}{l(P)!}
\bigwedge_{i=1}^{l(P)}\Omega_{p_i,p_i}\\&=& 
\sum_{P\in\nb_\star^{(\nb_\star)}\atop   |P|=d}  
\frac{1}{l(P)!\prod_{i=1}^{l(P)} p_i} \Omega^d=\Omega^d. 
\end{eqnarray*}
The last equality is proven noticing that the map $$
\begin{array}{ccc}
\Sigma_p&\to&\Sigma_p\\
\phi&\mapsto&\big(\phi(1),\phi(2),\cdots,\phi(p_1)\big)
\big(\phi(p_1+1),\phi(p_1+2),\cdots,\phi(p_2)\big) \cdots   
\end{array} $$
is surjective on the set of permutations  of type $P$ and is
$\prod_{i=1}^{l(P)} p_i$-to one. 
Each permutation of type $P$ is obtained by $l(P)!$ maps
of this kind.
\end{proof}

\bigskip
We can now compute the function $\Phi$ involved in the expression for
Bott-Chern forms. The first part of the proof of our main theorem will
be done in three steps. In the following proposition,
we first separate fiber differentials $dz$ and base differentials
arising in the curvature of $E^\star$
\begin{prop}
\label{prop:phi}
$$\Phi _{d+1}(u)=c'_{d}( E ^\star)+ \sum_{1\leq s\leq d} (1-u)^s 
\sum_{P,S\in\nb_\star^{(\nb_\star)}\atop  {l(P)=l(S), P\geq S\atop |P|\leq d,
  |S|=s, (l(P)\leq s)}}
\frac{ (-1)^{|P|+l(P)}}{l(P)!} c'_{d-|P|}(E
^\star)\bigwedge_{m=1}^{l(P)} \Omega_{p_m,s_m}.$$
\end{prop}

\begin{proof}
 First remark that 
\begin{eqnarray*}
\Phi _{d+1}(u)&=&\frac{1}{d!}\sum_{2\leq i_1,i_2,\cdots ,i_d\leq r}
\sum_{\sigma\in\Sigma_d}\ep (\sigma)\bigwedge _{m=1}^d 
(c_{i_m i_{\sigma (m)}} 
+ (1-u)\frac{i}{2\pi} dz_{i_m}\wedge d\zb_{i_{\sigma (m)}}).
\end{eqnarray*}
We first invert the summation over $\sigma\in\Sigma_d$ and the
summation resulting from the development of the above wedge
product. We then have to specify which $dz\wedge d\zb$ terms we
consider in the  development. The part $c'_{d}( E ^\star)$ is got when no
$dz\wedge d\zb$ terms are chosen. 
According to these, we write $\sigma$ as product of cycles 
with disjoint supports, neglecting the precise expression of those
cycles containing no $dz\wedge d\zb$ terms. If for example we seek for
the terms containing only $dz_{i_1}$ and $dz_{i_2}$ among the $dz$,
we will only write explicitly the cycles in $\sigma$ containing $1$ and $2$.
The permutation $\sigma$ is given by $P\in \nb_\star^{(\nb_\star)}$
with $|P|\leq d$ and $\phi :\coprod _{m=1}^{l(P)}
\nb_{p_m}^{(m)}\to\nb_d$ injective and $\sigma'$ permutation of
$\nb_d-im\phi$ by
$$\sigma=\sigma'\prod_{m=1}^{l(P)}
\left(\phi(1^{(m)}),\phi(2^{(m)}),\cdots,\phi(p_m^{(m)})\right).$$
Notice that each $\sigma$ is obtained $l(P)!\prod_{m=1}^{l(P)}p_m$
times. Hence,
$$\sum_{\sigma\in\Sigma_d}=\sum_{P} \frac{1}{l(P)!\prod_{m=1}^{l(P)}p_m}
\sum_{\phi}\sum_{\sigma'}.$$
Each $Q=\coprod _{m=1}^{l(P)}Q^{(m)}$, $Q^{(m)}\subset\nb_{p_m}$ gives the
summand ($\ep (\sigma)=\ep (\sigma ')(-1)^{|P|+l(P)}$)
$$(-1)^{|P|+l(P)}\bigwedge_{\tau =1\atop\tau\not\in im (\phi )}^d
\ep(\sigma') c_{i_\tau i_{\sigma'(\tau)}}\bigwedge_{m=1}^{l(P)}\bigwedge
_{t=1}^{p_m} C^{I_\varphi,Q^{(m)}}_{t^{(m)},(t+1)^{(m)}}.$$
 After commuting with
$\sum_{P}\sum_{\phi}\sum_{Q}$, rewrite
 $\displaystyle\sum_{2\leq i_1,i_2,\cdots ,i_d\leq r}$ as
$\displaystyle\sum_{2\leq i_\tau\leq r\atop \tau\not\in im(\varphi)}
\sum_{m=1}^{l(P)} 
\sum _{{{2\leq i_{\varphi(1^{(m)})}\leq r\atop 2\leq
  i_{\varphi(2^{(m)})}\leq r}\atop\vdots}\atop
  2\leq i_{\varphi(p_m^{(m)})}\leq r} $.  We get the expression for
$\Phi_{d+1}(u)$ as $c'_{d}( E ^\star)$ plus 
$$\frac{1}{d!}\sum_P\sum_{\phi}\sum_Q
\frac{(-1)^{|P|+l(P)}}{l(P)!\prod_{m=1}^{l(P)}p_m} 
(d-|P|)! c'_{d-|P|}(E^\star)
\bigwedge_{m=1}^{l(P)}(1-u)^{\sharp Q^{(m)}}\Omega_{p_m,Q^{(m)}}.$$
The summand now does not depend on $\phi$. Hence the summation over all
injective $\phi:\coprod _{m=1}^{l(P)}\nb_{p_m}^{(m)}\to\nb_d$ gives
the factor $\displaystyle\frac{d!}{(d-|P|)!}$. We infer
\begin{eqnarray}\label{formula}
\Phi_{d+1}(u)&=& c'_{d}( E ^\star)+\nonumber\\
&&\sum_P\sum_Q \frac{(-1)^{|P|+l(P)}}{l(P)!\prod_{m=1}^{l(P)}p_m}
c'_{d-|P|}(E^\star)\bigwedge_{m=1}^{l(P)}
(1-u)^{\sharp Q^{(m)}}\Omega_{p_m,Q^{(m)}}.
\end{eqnarray}
Summing now over all $Q$ with same size $S$ gives 
\begin{eqnarray*}
\Phi_{d+1}(u)&=&c'_{d}( E ^\star)+
\sum_P\sum_S \frac{(-1)^{|P|+l(P)}}{l(P)!\prod_{m=1}^{l(P)}p_m}
c'_{d-|P|}(E^\star) 
\bigwedge_{m=1}^{l(P)}(1-u)^{s_m} p_m\Omega_{p_m,s_m}\\
&=&c'_{d}( E ^\star)+
\sum_P\sum_S \frac{(-1)^{|P|+l(P)}}{l(P)!} c'_{d-|P|}(E^\star) (1-u)^{|S|}
\bigwedge_{m=1}^{l(P)} \Omega_{p_m,s_m}.
\end{eqnarray*}
\end{proof}

We now simplify the expression of fiber differentials.
Remark that when two $dz\wedge d\zb$ terms are neighbors ($q_{t+1}=q_t +1$), 
$$\Omega_{p,\{q_1,q_2,\cdots,q_{t},q_t +1,q_{t+2},\cdots ,q_s\}}=-\Omega 
\Omega_{p-1,\{q_1,q_2,\cdots,q_{t-1},q_t,q_{t+2}-1,\cdots ,q_s-1\}}
$$
This remark will be improved in the following proposition.
\begin{prop}
For all $(p,s)\in\nb^2$ with $s\leq p$,
$$\Omega_{p,s}=\sum_{1\leq b_1\leq b_2\leq\cdots\leq b_s\leq
  p \atop |B|=p}\frac{(-1)^{s-1}(s-1)!}{h(B)!}
 \Omega_{b_1,1}\Omega_{b_2,1}\cdots \Omega_{b_s,1}.$$
\end{prop}

\begin{proof}
Recall that 
$\Omega_{p,s}:=\frac{1}{p}\sum_{Q\subset\nb_p\atop \sharp Q=s}
\sum_{2\leq i_1,i_2,\cdots,i_p\leq r}\bigwedge_{a=1}^p
C_{a,a+1}^{I,Q}$.
Gathering terms in the sum in the following way
\begin{eqnarray*}
\big(d\zb_{i_{q_1+1}}c_{i_{q_1+1} i_{q_1+2}}&\cdots & c_{i_{q_2-1}i_{q_2}} 
dz_{i_{q_2}}\big)\cdots
\big(d\zb_{i_{q_{s-1}+1}}c_{i_{q_{s-1}+1} i_{q_{s-1}+2}}\cdots 
c_{i_{q_s-1}i_{q_s}} dz_{i_{q_s}}\big) \\ &&
\big(-d\zb_{i_{q_{s}+1}}c_{i_{q_{s}+1} i_{q_{s}+2}}\cdots
c_{i_{p-1}i_p}c_{i_p i_1}c_{i_1i_2}\cdots
c_{i_{q_1-1}i_{q_1}}dz_{i_{q_1}}\big) 
\end{eqnarray*}
we infer that each $Q$ contributes to
$\big(-\Omega_{q_2-q_1,1}\big)\big(-\Omega_{q_3-q_2,1}\big)\cdots 
\big(-\Omega_{q_s-q_{s-1},1}\big)
\big(~+~\Omega_{p-q_s+q_1,1}\big)$.
Set $b'_j:=q_{j+1}-q_j$ for $1\leq j\leq s-1$ and $b'_s:=p-q_s+q_1$,
we get a map ($\mathcal{P}_s(\nb_p)$ denotes the set of subsets of
$\nb_p$  of cardinal $s$)
$$
\begin{array}{ccccc}
\mathcal{P}_s(\nb_p)&\stackrel{\alpha}\to&(\nb_\star)^s 
&\stackrel{\beta}\to &
\{\textrm{ partitions of weight } p \textrm{ and length } s\}\\
Q&\mapsto &(b'_1,b'_2,\cdots ,b'_s)&\mapsto &b_1\leq b_2\leq \cdots
\leq b_s
\end{array}$$
The fiber $\alpha ^{-1}(B')$ has $b'_s$ elements and the fiber $\beta
^{-1}(B)$ $\frac{s!}{h(B)!}$. But
$$\sum_{B'\in\beta^{-1}(B)}b'_s=\frac{1}{s}\sum_{B'\in\beta^{-1}(B)}
b'_1+b'_2+\cdots +b'_s=\frac{p}{s}\frac{s!}{h(B)!}.$$ 
Hence the composed map is of degree $\displaystyle\frac{p
  (s-1)!}{h(B)!}$ over $B$.
\end{proof}

\bigskip
We are now about the final step
\begin{prop}
\label{prop:fin} For all $d\in \nb$,
$$\Phi _{d+1}(u)=c'_{d}( E ^\star)+ \sum_{1\leq s\leq p \leq
  d}\sum_{1\leq b_1\leq b_2\cdots \leq b_s\leq p\atop |B|=p}
(1-u)^s \frac{(-1)^{p+s}s!}{h(B)!}
c'_{d-p}(E ^\star)\Omega_{b_1,1}\Omega_{b_2,1}\cdots \Omega_{b_s,1}$$
\end{prop}

\begin{proof}
Fix $s\leq p$ in $\nb_\star$.
The map  which assigns to the product $Q$ of
$Q^{(m)}\in\mathcal{P}_{s_m} (\nb_{p_m})$ the associated product of
partitions $B^{(m)}$ is $\displaystyle\prod_{m=1}^{l(P)}
\frac{p_m (s_m-1)!}{h(B^{(m)})!}$ to $1$. 
The map  which assigns to
  the latter product of partitions a partition $B$ of weight $p$ and of
  length $s$ by concatenation and reordering is
  $\displaystyle\sum_{S\in(\nb_\star)^{(\nb_\star)}\atop |S|=s  }\frac{s!
    \prod_{m=1}^{l(S)}
    h(B^{(m)})!}{h(B)!\prod_{m=1}^{l(S)}s_m!}$ to $1$. Hence, the composed
  map is $\displaystyle\sum_{S\in(\nb_\star)^{(\nb_\star)}\atop |S|=s
  }\frac{s! 
    \prod_{m=1}^{l(S)}p_m }{h(B)!\prod_{m=1}^{l(S)}s_m}$ to $1$.
Back to the formula~(\ref{formula}), we obtain a new expression for 
$\Phi_{d+1}(u)$.
\begin{eqnarray*}
\Phi_{d+1}(u)&=&c'_{d}( E ^\star)+
\sum_P\sum_Q \frac{(-1)^{|P|+l(P)}}{l(P)!\prod_{m=1}^{l(P)}p_m}
c'_{d-|P|}(E^\star)\bigwedge_{m=1}^{l(P)}
(1-u)^{\sharp Q^{(m)}}\Omega_{p_m,Q^{(m)}}\\
&=& c'_{d}( E ^\star)
+\sum_{1\leq s\leq p \leq
  d}\sum_{1\leq b_1\leq b_2\cdots \leq b_s\leq p\atop |B|=p}
\sum_{S\in(\nb_\star)^{(\nb_\star)}\atop |S|=s}c'_{d-p}(E^\star)
(1-u)^s\\&&
\frac{s!\prod_{m=1}^{l(S)}p_m }{h(B)!\prod_{m=1}^{l(S)}s_m}
\frac{(-1)^{p+l(S)}}{l(S)!\prod_{m=1}^{l(S)}p_m}
 (-1)^{s+l(S)}
\Omega_{b_1,1}\Omega_{b_2,1}\cdots \Omega_{b_s,1}.
\end{eqnarray*} 
We recall once more the identity
$\displaystyle\sum_{S\in(\nb_\star)^{(\nb_\star)}\atop |S|=s}
\frac{1}{l(S)!\prod_{m=1}^{l(S)}s_m} =1$ to end the proof.
\end{proof}

\bigskip
Hence, we get the Bott-Chern forms by integrating.
\begin{theo}\label{theo}
\begin{eqnarray*}
\widetilde{c_{d+1}}=- 
 \sum_{1\leq s\leq p \leq d}
\sum_{1\leq b_1\leq b_2\cdots \leq b_s\leq p\atop |B|=p}
\mathcal{H}_s  \frac{(-1)^{p+s}s!}{h(B)!}
c'_{d-p}(E ^\star)\Omega_{b_1,1}\Omega_{b_2,1}\cdots \Omega_{b_s,1}.
\end{eqnarray*}
\end{theo}
Here $\mathcal{H}_s$ is the harmonic number $\displaystyle\sum_1^s
 \frac{1}{i}=\int_0^1 \frac{1-(1-u)^s}{u}$.
This formula reduces in the case of flat vector bundle to that 
of Gillet and Soul{\'e} (proposition 5.3 of ~\cite{gs}).
Coordinates free expressions for $c'_d(E ^\star)$ and $\Omega_{b,1}$
will be given in proposition~\ref{prop:c} and ~\ref{prop:omegas}.
Notice that the relative degree of each summand is $2s$ and 
the degree in base variables is 
 $2(d-p+|B|-s)=2(d-s)$. We can therefore restrict the range of $s$ 
to $\max(d-dim X,1)\leq s\leq \min (r-1,d)$.
\bigskip

\bigskip\bigskip
\section{Coordinates free expressions}
\label{sec:coordinates}

We now give coordinates free expressions of $c'_i(E ^\star)$
and of the $\displaystyle\Omega_{p,1}$. We will express the results
in terms of the following forms in $A^{q,q}(\pb (E),\cb)$~: for
$q\in\nb $,
 at the point $(x,[a^\star ])\in \pb (E)$,
$$\Theta ^q a^\star:=
\frac{\langle\pi^\star\Theta(E ^\star,h)^q a ^\star,a
 ^\star\rangle}{\n a ^\star\n^2}.$$
Recall that $\Theta(E ^\star,h)$ is in $A^{1,1}(X,End(E))$.
The $q$-th power is taken with the wedge product in the form part and
the composition in the endomorphism part. We will also need
the form $\alpha=\Theta(\oc_E(1),h)$.

\bigskip
\begin{prop}\label{prop:c}
$$c'_d(E ^\star)=\sum_{p+m=d\atop p,m\geq 0}(-1)^{p}(\Theta^p
a^\star)\pi^\star c_m(E ^\star,h) $$
where we have set $(\Theta^0 a^\star)=\pi^\star c_0(E ^\star)=1$. 
\end{prop}

\begin{proof}
The expression for $\Theta ^q a^\star$ in coordinates is 
$$\frac{ \sum_{J\in\{1,\cdots,r\}^{q+2}}
\bigwedge_{m=1}^{q} c_{j_m j_{m+1}}a_{j_{q+1}}<e ^\star_{j_{1}}, e
^\star _{j_{q+2}}> \overline{a_{j_{q+2}}}}
{\sum_{(j,j')\in\{1,\cdots,r\}^{2}} a_{j}<e ^\star_{j}, e
^\star _{j'}> \overline{a_{j'}} }$$
and  is $\sum_{J\in\{1,\cdots,r\}^{q+1}\atop j_1=j_{q+1}=1}
\prod_{a=1}^{q} c_{j_a j_{a+1}}$ at the center $x_0$ of the
normal frame for $(E,h)$. Hence,
\begin{eqnarray*}
c'_1(E ^\star)&=&
\pi^\star c_1(E ^\star,h)-(\Theta ^1 a^\star) \\
c'_2(E ^\star)&=&\pi^\star c_2(E ^\star,h)-\left[c_{11}c'_1(E
  ^\star)-\sum_{j\geq 2}c_{1j}c_{j1}\right]\\
&=&\pi^\star c_2(E ^\star,h)-(\Theta ^1 a^\star)\pi^\star c_1(E
^\star,h)+(\Theta ^2 a^\star).
    \end{eqnarray*}
In the same spirit of the proof of proposition~\ref{prop:phi}, we
write the permutations according to one of the positions of $1$ (if any)~: 
$\sigma=(\phi(1),\phi(2),\cdots,\phi(p))\sigma'$ if $i_{\phi (1)}=1$.
\begin{eqnarray*}
\lefteqn{c_d(E ^\star,h)= c'_d(E ^\star)}\\&&+\frac{1}{d!}
\sum_{1\leq i_1,i_2,\cdots ,i_d\leq r\atop 1\in\{i_1,i_2,\cdots ,i_d\}}
\sum_{p=1}^d\sum_{{\phi:\nb_p\to\nb_d\atop \textrm{injective}}\atop i_{\phi
      (1)}=1}\sum_{\sigma '\in\Sigma_{d-p}}
 (-1)^{p+1} c_{1i_{\phi(2)}}
c_{i_{\phi(2)} i_{\phi(3)}}\cdots
c_{i_{\phi(p)}1}\!\!\bigwedge_{\tau\not\in im\phi}
\!\!c_{i_\tau,i_{\sigma'(\tau)}}. 
 \end{eqnarray*}
Rewriting $\sum\limits_{1\leq i_1,i_2,\cdots ,i_d\leq r\atop
  1\in\{i_1,i_2,\cdots,i_d\} }$ as $\sum\limits_{1\leq
  i_{\phi(2)},\cdots ,i_{\phi(p)}\leq r}\sum\limits_{1\leq i_\tau\leq r\atop 
\tau\not\in im\phi}$, then replacing the sum
  $\sum\limits_{\phi:\nb_p\to\nb_d\textrm{injective},\atop i_{\phi (1)}=1}$
  by product by factor $\displaystyle\frac{d!}{d-p)!}$ we infer
\begin{eqnarray*}
{c_d(E ^\star,h)}&=&c'_d(E ^\star)+\sum_{p=1}^d 
(-1)^{p+1} (\Theta^p a^\star) \pi^\star c_{d-p}(E ^\star,h).
\end{eqnarray*}
\end{proof}

\bigskip
Before computing the $\Omega_{p,1}$, we need a lemma which relies
on special properties of the normal frame.
\begin{lemma}\label{normal}
In the center of a normal frame for an holomorphic  Hermitian vector
bundle $(E,h)$, $d\Theta _E=0$ and $id'd''\Theta _E=0$.
\end{lemma}
\begin{proof}
Call $H$ the metric matrix and $A$ the connection $1$-form for the
Chern connection of $(E,h)$ 
in a frame normal at $x_0$. By the compatibility of the connection
with the metric and the holomorphic structure, we infer that
$d'H=AH$. Then $A(x_0)=0$. On the whole chart, 
$\Theta _E = \frac{i}{2\pi} (dA+A\wedge A)$. So, $d\Theta _E =\Theta
_E \wedge A -A\wedge\Theta _E$ (Bianchi identity). Hence, $d\Theta
_E(x_0)=0$. Now, at $x_0$
\begin{eqnarray*}
d'd''\Theta _E(x_0)=d'd\Theta _E(x_0)&=&d'\Theta_E \wedge A +\Theta_E \wedge
d'A-d'A\wedge\Theta _E+ A\wedge d'\Theta _E\\
&=&  \Theta_E \wedge d'A-d'A\wedge\Theta _E
\end{eqnarray*}
From, $d'H=AH$ and $A(x_0)=0$ we get $d'A(x_0)=0$ (Notice that
$\Theta_E(x_0)= dA(x_0)$ is of type $(1,1)$).
\end{proof}

\bigskip
\begin{prop}
\label{prop:omegas} For $d\geq 1$, at the point $(x_0,[a^\star_0])$,
\begin{eqnarray*}
\Omega_{d+1,1}&=&\frac{i}{2\pi}d'd''(\Theta ^d a^\star) 
+ (\Theta ^{d+1} a^\star)
+ \alpha (\Theta ^d a^\star) \\
&+&\frac{i}{2\pi}\sum_{m=2}^d \sum_{B\in{\nb_\star}^{(\nb_\star)}\atop
  l(B)=m,|B|=d}  
(-1)^m (d''\Theta ^{b_1}a ^\star)(\Theta ^{b_2}a ^\star)
(\Theta^{b_3}a ^\star)\cdots (\Theta ^{b_{m-1}}a ^\star)
(d'\Theta ^{b_m}a ^\star)\end{eqnarray*}
\end{prop}

\begin{proof}
From the definition of the normal frame we get at $x_0$, $d\n a
^\star\n^2 =0$. Hence, at $(x_0,[a_0^\star])$, 
\begin{eqnarray*}
\frac{i}{2\pi}d'd''(\Theta ^d a^\star) &=&
\frac{i}{2\pi}d'd''\langle\pi^\star\Theta(E ^\star,h)^d a ^\star,a^\star\rangle
-(\Theta ^d a^\star)\alpha.
\end{eqnarray*}
From lemma~\ref{normal} and the relation
$\frac{i}{2\pi}d'd''{\langle e ^\star_i,e^\star_j\rangle}
=-c_{ji}$ we infer using the expression of 
$\langle\pi^\star\Theta(E ^\star,h)^d a^\star,a^\star\rangle $  in
coordinates  
\begin{eqnarray*}
\lefteqn{\frac{i}{2\pi}d'd''\langle\pi^\star\Theta(E ^\star,h)^d
 a^\star,a^\star\rangle }\\&=&\frac{i}{2\pi}
 \sum_{J\in\{1,\cdots,r\}^{d+1}}
\bigwedge_{m=1}^{d} c_{j_m j_{m+1}}dz_{j_{d+1}}\wedge
d\overline{z_{j_{1}}} 
+ \sum_{J\in\{1,\cdots,r\}^{d+2}\atop j_{d+1}=j_{d+2}=1 }
\bigwedge_{m=1}^{d} c_{j_m j_{m+1}}(-c_{j_{d+2} j_{1} })\\
&=&\frac{i}{2\pi}
 \sum_{J\in\{1,\cdots,r\}^{d+1}}
\bigwedge_{m=1}^{d} c_{j_m j_{m+1}}dz_{j_{d+1}}\wedge
d\overline{z_{j_{1}}} -(\Theta ^{d+1} a^\star).\\
\end{eqnarray*}
Now, to recover $\Omega_{d+1,1}$ we have to remove the value $1$ for
the $j$'s in the first sum. Gathering the terms according to the fact
that they contain at least $l$ $j's$ achieving the 
 value $1$ ($j_{d+1}\not =1,j_{1}\not =1$), we get (notice that we
 will interchange $dz_{j_{d+1}}$ and~$d\overline{z_{j_1}} $)
\begin{eqnarray*}
\lefteqn{
\frac{i}{2\pi}\sum_{J\in\{1,\cdots,r\}^{d+1}}
\left(\bigwedge_{m=1}^{d} c_{j_m j_{m+1}}\right)
dz_{j_{d+1}}\wedge d\overline{z_{j_1}}  }\\
&=&\Omega_{d+1,1}+\frac{i}{2\pi}\sum_{l=1}^{d-1}(-1)^l
\sum_{1<q_1<q_2<\atop \cdots < q_l<d+1}
\sum_{J-Q\in\{1,\cdots,r\}^{d+1-l}}
(d\overline{z}_{j_1} c_{j_1j_2}c_{j_2j_3}\cdots  c_{j_{q_1-1}1})\\&&
(c_{1j_{q_1+1}}\cdots c_{j_{q_2-1}1})(c_{1j_{q_2+1}}\cdots
c_{j_{q_3-1}1})
\cdots
(c_{1j_{q_{l-1}+1}}\cdots c_{j_{q_l-1}1})
(c_{1j_{q_l+1}}\cdots c_{j_dj_{d+1}}dz_{j_{d+1}})\\
&=&\Omega_{d+1,1}+\frac{i}{2\pi}\sum_{l=1}^{d-1}(-1)^l
 (\sum_{1\leq j_1,j_2,\cdots , j_{q_1-1}\leq r}
  d\overline{z}_{j_1} c_{j_1j_2}c_{j_2j_3}\cdots  c_{j_{q_1-1}1}) \\
&&(\Theta ^{q_2-q_1} a^\star)(\Theta ^{q_3-q_2} a^\star)
\cdots (\Theta ^{q_{l}-q_{l-1}} a^\star)
(\sum_{1\leq j_{q_l+1},\cdots , j_{d+1}\leq r}
 c_{1j_{q_l+1}}\cdots c_{j_dj_{d+1}}dz_{j_{d+1}})
\end{eqnarray*}
Using once more the normality of the frame, we compute at the point
$(x_0,[a_0])$ 
\begin{eqnarray*}
\sum_{1\leq j_1,j_2,\cdots , j_{q_1-1}\leq r}
  d\overline{z}_{j_1} c_{j_1j_2}c_{j_2j_3}\cdots  c_{j_{q_1-1}1}
&=&d''(\Theta ^{q_1-1} a^\star) \\
\sum_{1\leq j_{q_l+1},\cdots , j_{d+1}\leq r}
 c_{1j_{q_l+1}}\cdots c_{j_dj_{d+1}}dz_{j_{d+1}}
&=& d'(\Theta ^{d+1-q_l} a^\star).
\end{eqnarray*}
Now, the proof will be complete if we set $b_1:=q_1-1\in\nb_\star$,
$b_j:=q_j-q_{j-1}$ for $j$ between~$2$ and $l$ and
$b_{l+1}:=d+1-q_l\in\nb_\star$ and $m=l+1$.
\end{proof}

\bigskip\bigskip
\section{On the Bott-Chern forms generating polynomial}
Bott-Chern forms fulfill a recursion formula which is easily expressed
in terms of their generating polynomial. The latter has a concise
expression.

Set 
 \begin{eqnarray*}
c_t(E ^\star,h):=\sum_{d=0}^{+\infty} t^d c_d(E ^\star,h) \ ;\
\widetilde{c_t}(\Sigma,h):=\sum_{d=0}^{+\infty} t^d
\widetilde{c_{d+1}}(\Sigma,h) \  ;\   
\Theta _t:=\sum_{d=0}^{+\infty}(-t)^d (\Theta ^d a ^\star).
\end{eqnarray*}
 We will also need as auxiliary notations
 \begin{eqnarray*}
c'_t(E^\star):=\sum_{d=0}^{+\infty} t^d c'_d(E ^\star)\ ;\
\Omega_t:= \sum_{d=1}^{+\infty}(-t) ^d\Omega _{d,1} \ ;\
H(X):=\sum_{s=1}^{+\infty}\mathcal{H}_s X^{s}.
\end{eqnarray*}
Theorem~\ref{theo} can be rewritten as
\begin{eqnarray*}
t^d\widetilde{c_{d+1}}=- 
 \sum_{1\leq s \leq d}
\sum_{B\in(\nb_\star)^s\atop e+\sum b_i=d}
\mathcal{H}_s  {(-1)^{s}}
t^{e} c'_{e}(E ^\star)
(-t)^{b_1}\Omega_{b_1,1}(-t)^{b_2}\Omega_{b_2,1}\cdots 
(-t)^{b_s}\Omega_{b_s,1}.
\end{eqnarray*}
Summing over $d$, we infer
\begin{eqnarray*}
\widetilde{c_t}(\Sigma,h)=-\sum_{1\leq s}\mathcal{H}_s  {(-1)^{s}}
c'_t(E^\star) (\Omega_t)^s=-H(-\Omega_t)c'_t(E^\star).
\end{eqnarray*}
Proposition~\ref{prop:c} can be rewritten as
$c'_t(E^\star)=\Theta _t\pi ^\star c_t(E ^\star,h)$.
As for proposition~\ref{prop:omegas} notice that
\begin{eqnarray*}
\sum_{d=0}^{+\infty}(-t)^{d+1}
\frac{i}{2\pi}\sum_{m=2}^d \sum_{B\in{\nb_\star}^{(\nb_\star)}\atop
  l(B)=m,|B|=d}  
(-1)^m (d''\Theta ^{b_1}a ^\star)(\Theta ^{b_2}a ^\star)
(\Theta^{b_3}a ^\star)\cdots (\Theta ^{b_{m-1}}a ^\star)
(d'\Theta ^{b_m}a ^\star)\\
=-t\frac{i}{2\pi}(\sum_{b_1=1}^{+\infty}(-t)^{b_1}d''\Theta ^{b_1}a ^\star)
(\sum_{b_m=1}^{+\infty}(-t)^{b_m}d'\Theta ^{b_m}a ^\star)
\sum_{m=2}^{+\infty}(1-\Theta_t)^{m-2}\\
=-t\frac{i}{2\pi}\frac{d''\Theta_t d'\Theta_t}{\Theta_t}
\end{eqnarray*}
so that
\begin{eqnarray*}
-\Omega_t=t\left(\frac{i}{2\pi}d'd''\Theta_t -\frac{\Theta_t-1}{t}
+\alpha\Theta_t+\frac{i}{2\pi}\frac{d''\Theta_t d'\Theta_t}{\Theta_t}
\right).
\end{eqnarray*}
We are led to our main theorem
\begin{theo}\label{theo:total} Bott-Chern forms for the metric Euler
  sequence can be chosen to be
$$\widetilde{c_t}(\Sigma,h)=-H\left(t\alpha\Theta_t
+t\frac{1-\Theta_t}{t}
+t\Theta_t\frac{i}{2\pi}d'd''\ln\Theta_t \right)
\Theta_t\pi ^\star c_t(E ^\star,h).$$
\end{theo}

\bigskip\bigskip
\section{Some special cases}

\bigskip
\subsection{On first parts of Bott-Chern forms}
The high relative degree part of Bott-Chern forms is easy to compute
with our formulas. Half of the relative degree will be
indicated by an extra indices. Theorem~\ref{theo} reads
 \begin{eqnarray*}
(\widetilde{c_{d+1}})_{d}&=&-\mathcal{H}_d \Omega ^d \\
(\widetilde{c_{d+1}})_{d-1}&=&-\mathcal{H}_{d-1}\left( 
c'_1(E ^\star)\Omega ^{d-1}
- (d-1) \Omega_{2,1}\Omega ^{d-2}\right)\\
(\widetilde{c_{d+1}})_{d-2}&=&-\mathcal{H}_{d-2}\Big(
c'_2(E ^\star)\Omega ^{d-2}
- (d-2) c'_1(E ^\star)\Omega_{2,1}\Omega ^{d-2}\\
&&
+(d-2) \Omega_{3,1}\Omega ^{d-3}
+\frac{(d-2)(d-3)}{2}(\Omega_{2,1})^2\Omega ^{d-4}\Big).
\end{eqnarray*}

\bigskip
\subsection{On first Bott-Chern forms}

The propositions~\ref{prop:c} and \ref{prop:omegas}
give theorem~\ref{theo} an intrinsic form. Alternatively, we may use
theorem~\ref{theo:total}.
Explicitly, we get up to  $d'$-exact and $d''$-exact forms, 
\begin{eqnarray*}
\widetilde{c_1}&=& 0\\
\widetilde{c_2}&=& -\left[\mathcal{H}_1 \alpha\right] -(\Theta ^1 a^\star)\\
\widetilde{c_3}&=& -\left[\mathcal{H}_2\alpha ^2 +\pi ^\star c_1(E
^\star,h)\alpha\right] 
  -\left[\alpha+\pi ^\star c_1(E^\star,h)\right]  (\Theta^1 a^\star) 
-\left(\frac{1}{2} (\Theta^1 a^\star)^2- (\Theta^2 a^\star)\right).
\end{eqnarray*}
The terms in bracket are $d'd''$-closed. 

\bigskip
\subsection{For curves and surfaces}
\label{sec:curve}

We now assume that the manifold $X$ is a curve or a surface. Only the
high relative degree part will occur in theorem~\ref{theo}.
Now, make use of
\begin{eqnarray*}
 \Omega&=& \alpha +(\Theta ^1 a^\star)\\
\Omega_{2,1}&=&\frac{i}{2\pi}d'd''(\Theta ^1 a^\star) 
+(\Theta ^2 a^\star)+\alpha(\Theta ^1 a^\star)\\
\Omega_{3,1}&=& \frac{i}{2\pi}d'd''(\Theta ^2 a^\star)
+(\Theta ^3 a^\star)+\alpha(\Theta ^2 a^\star)
+\frac{i}{2\pi}d''(\Theta ^1 a^\star) d'(\Theta ^1 a^\star).
\end{eqnarray*}

Hence we proved that
 up to $d'$-exact and $d''$-exact forms Bott-Chern forms are given on
 curves by
\begin{eqnarray*}
\widetilde{c_{d+1}}&=&-\left[\mathcal{H}_d \alpha ^d 
+\mathcal{H}_{d-1}\pi ^\star c_1(E ^\star ,h)\alpha ^{d-1}\right]
-\left[\alpha ^{d-1}\right](\Theta ^1 a^\star)
\end{eqnarray*}
and on surfaces by
\begin{eqnarray*}
\widetilde{c_{d+1}}&=&-\left[\mathcal{H}_d \alpha ^d 
+\mathcal{H}_{d-1}\pi ^\star c_1(E ^\star ,h)\alpha ^{d-1}
+\mathcal{H}_{d-2}\pi ^\star c_2(E ^\star ,h)\alpha ^{d-2}\right]\\
&&-\left[ \alpha ^{d-1}+\pi ^\star c_1(E ^\star ,h)\alpha
  ^{d-2}-(d-2)\alpha ^{d-3}\frac{i}{2\pi}d'd''\Theta ^1 a^\star\right]
(\Theta ^1 a^\star)\\
&&-\left[\alpha ^{d-2}\right]\left( \frac{1}{2} (\Theta^1
  a^\star)^2-(\Theta^2 a^\star)\right).
\end{eqnarray*}

\bigskip\bigskip
\section{On  characteristic forms}
\label{form}

\bigskip
\subsection{From Segre forms to Chern forms}

The geometric Segre forms are defined for $i~\in~\nb$ by
 $$s'_i(E,h):=\pi_\star (\alpha ^{r-1+i})
\in A^{i,i} (X,\cb).$$
Consider the metric Euler exact 
sequence~$(\Sigma,h)$ twisted by $\oc _E(1)$ 
\begin{eqnarray*}
  0\to \oc_{\pb (E)}\to\pi^\star E^\star\otimes\oc _E(1)\to T_{\pb (E)/X}\to 0 
\end{eqnarray*}
and its top Bott-Chern form $\widetilde{c_{r}}(\Sigma(1),h)$.
Notice that $c_t(\oc_{\pb (E)})=1$ for the induced metric on
$\oc_{\pb (E)}$ is the flat one. Degree considerations lead to the
 relation 
$$c_r(\pi ^\star E ^\star\ot\oc _E(1),h) 
=-\frac{i}{2\pi}d'd''\left(\widetilde{c_{r}}(\Sigma(1),h)\right).$$
which in cohomology reduces to Grothendieck defining relation 
for Chern classes of~$E$.
Now, for $m\in\nb$ set $S_{m+1} (E,h):=\pi_\star\left(\alpha^{m}
\widetilde{c_{r}}(\Sigma(1),h) \right)\in \widetilde{A^{m,m}}(X).$
Define the class $R_{m+1}(E,h)\in \widetilde{A^{m,m}}(X)$ by
$$\sum_{m=0}^{+\infty}t^m R_{m+1}(E,h)=
\left(\sum_{m=0}^{+\infty} (-t)^m c_m(E,h) \right)^{-1}
\left(\sum_{m=0}^{+\infty}t^m S_{m+1} (E,h)\right).$$

\bigskip
For $m\geq 1$, using the projection formula, we infer
\begin{eqnarray*}
  \lefteqn{\sum_{p+q=m \atop p,q\in\nb} c_p( E ^\star ,h)s'_q(E,h)}
\\
&=& 
\pi_\star \left(\sum_{p+q=m \atop p,q\in\nb}  
c_p(\pi^\star E^\star,h) \alpha ^{r-1+q}\right)
= \pi_\star \left(\sum_{p+q=m \atop p,q\in\zbb}  
c_p(\pi^\star E^\star,h) \alpha ^{r-1+q}\right)\\
&=&\pi_\star\left( \alpha ^{m-1}\sum_{p+s=r \atop p,s\in\zbb}  
c_p(\pi^\star E^\star,h) \alpha ^{s}\right)
=\pi_\star\left( \alpha ^{m-1} c_r(\pi^\star E^\star\otimes\oc_E(1)
  ,h)\right)\\
&=& -\frac{i}{2\pi} d'd''\pi_\star\left(\alpha^{m-1}
\widetilde{c_{r}}(\Sigma(1),h) \right).
\end{eqnarray*} 
From the previous equation, we infer that the Chern forms of $(E,h)$
are related to the geometric Segre forms $s'$ by 
$$\left(\sum_{p=0}^{+\infty} (-t)^p c_p(E,h) \right)
\left(\sum_{q=0}^{+\infty} t^q s'_q(E,h)\right)= 
1-\frac{it}{2\pi} d'd''\sum_{m=0}^{+\infty}t^m S_{m+1} (E,h)
$$
that is
\begin{eqnarray}\label{arithmetique}
\\ \nonumber
\left( \sum_{m=0}^{+\infty} (-t)^m c_m(E,h) \right)^{-1}=
1+\sum_{m=1}^{+\infty} t^m \left( s'_m(E,h)
+\frac{i}{2\pi} d'd'' R_{m}(E,h)\right).
\end{eqnarray} 

\bigskip
The striking fact is that despite the appearance of the non-closed
forms $\Theta ^q a ^\star$ in the expression of
$\widetilde{c_{r}}(\Sigma(1),h)$  the forms $R$ are $d'd''$-closed.

\begin{prop}
\label{prop:closed}
The forms $R$ are $d'd''$-closed or equivalently the Chern forms are
related to the geometric Segre forms $s'$ by
$$\left(\sum_{m=0}^{+\infty} c_m(E^\star,h)\right)
\left(\sum_{m=0}^{+\infty} t^m s'_m(E,h)\right)=1$$
\end{prop}

\begin{proof}
We first compute the Segre forms.
\begin{eqnarray*}
\lefteqn{s'_m(E,h)=\pi_\star(\alpha ^{r-1+m})=\left({r-1+m\atop m}\right)
\pi_\star \left(\Omega ^{r-1}(-\Theta ^1a ^\star)^m\right)}\\
&=&\left({r-1+m\atop m}\right)(-1)^m \!\! \int\limits_{\pb^{r-1}}\Omega
^{r-1}\frac{ (\sum c_{ij}a_j \overline{a_i})^m}{|a|^{2m}}\\
&=&\left({r-1+m\atop m}\right)(-1)^m \!\!\!\!\!\!\sum_{I,J\in(\nb_r)^m}
\!\!\!\! c_{i_1j_1} c_{i_2j_2} \cdots c_{i_mj_m}
\!\!\int\limits_{\pb^{r-1}}\Omega^{r-1}
\frac{a_{j_1}\overline{a_{i_1}}a_{j_2}\overline{a_{i_2}}\cdots
a_{j_m}\overline{a_{i_m}}}{|a|^{2m}} 
\end{eqnarray*}
For parity reasons, the integral of non real terms vanishes. We may
hence restrict to $I=J$ as sets with multiplicities. 
That is there exists a permutation
$\sigma\in\Sigma_m$ such that for each
$\lambda$, $j_\lambda=i_{\sigma(\lambda )}$.
 One term may be gotten from different permutations
if some $i_\lambda$ equals some $i_\mu$ for $\lambda\not =\mu$. As in
the case of partitions,
we define the height $h(I)$  of a sequence of numbers $I$ to be the sequence
of cardinals of subsets of identical entries.
The term $a_{j_1}\overline{a_{i_1}}a_{j_2}\overline{a_{i_2}}\cdots
a_{j_m}\overline{a_{i_m}}$ can be written as 
$a_{i_{\sigma(1)}}\overline{a_{i_1}}a_{i_{\sigma(2)}}\overline{a_{i_2}}\cdots
a_{i_{\sigma(m)}}\overline{a_{i_m}}$ for $h(I)!$-different permutations.
We will make use of the formula
\begin{eqnarray*}\label{integral}
\int_{\pb^{r-1}} \frac{|a_1|^{2m_1}|a_2|^{2m_2}\cdots
  |a_r|^{2m_r}}{|a|^{2m}} \Omega ^{r-1}
=\frac{(r-1)!\prod m_i!}{(r-1+m)!}
\end{eqnarray*}
 where $m_i\in\nb$ are such that $\sum m_i=m$.
Back to our computations, we get
  \begin{eqnarray*}
s'_m(E,h)
&=&\left({r-1+m\atop m}\right)(-1)^m
\sum_{I\in(\nb_r)^m }\sum_{\sigma\in\Sigma_m}\\
 &&c_{i_1i_{\sigma(1)}}c_{i_2i_{\sigma(2)}} \cdots c_{i_mi_{\sigma(m)}}
\int_{\pb^{r-1}}\frac{\Omega^{r-1}}{h(I)!}\frac{|a_{i_1}|^2|a_{i_2}|^2\cdots
  |a_{i_m}|^2} {|a|^{2m}}
 \\
&=&\frac{(-1)^m}{m!}
\sum_{I\in(\nb_r)^m}\sum_{\sigma\in\Sigma_m} 
 c_{i_1i_{\sigma(1)}}c_{i_2i_{\sigma(2)}} \cdots c_{i_mi_{\sigma(m)}}.
\end{eqnarray*}

For a finite sequence $N$ in $\nb ^{(\nb_\star)}$ of non negative integers,
a permutation $\sigma$ in $\Sigma_{m}$ will be said to be of
\textit{shape $1^{n_1} 2^{n_2} \cdots m^{n_m}$},  ($\sum p n_p=m$) 
if it can be written as the product
of $n_p$ cycles of length $p$ with disjoint support.
There are $\frac{m!}{ \prod_p n_p! p^{n_p} }$ permutations of shape
$1^{n_1} 2^{n_2} \cdots m^{n_m}$ in $\Sigma_{m}$ and each will give the same
contribution in the sum after having computed $\sum_{I\in(\nb_\star)^m }$.
For $p\in\nb$, consider the closed forms on $X$ given by 
$$\theta_p:=Trace(\Theta(E,h)^{\ot p})
=(-1)^p \sum_{I\in(\nb_r)^p}c_{i_1i_2}c_{i_2i_3}\cdots 
c_{i_pi_1}.$$  
We now consider the total Segre form.
\begin{eqnarray*}
\sum_{m\in\nb} t^m s'_m(E,h)&=&\sum_m t^m\frac{1}{m!}
 \sum_{N\in\nb^{dim X}\atop \sum pn_p=m}
\prod _{p=1}^{dim X} \theta _p^{n_p} \frac{m!}{\prod_p n_p!p^{n_p}}\\
&=&\sum_{N\in\nb^{dim X}}\prod _{p=1}^{dim X}\left(\frac{t^p
 \theta_p}{p}\right)^{n_p}\frac{1}{n_p!}
=\prod _{p=1}^{dim X}\exp \left(\frac{t^p \theta_p}{p}\right).
\end{eqnarray*}
Noticing that the signature of a permutation of shape  $1^{n_1}
2^{n_2} \cdots m^{n_m}$ is $(-1)^{\sum_p n_p (p+1)}$, we also get 
\begin{eqnarray*}
\sum_{m\in\nb} t^m c_m(E ^\star,h)&=&
\sum_{N\in\nb^{dim X}}\prod _{p=1}^{dim X}\left(-\frac{t^p
 \theta_p}{p}\right)^{n_p}\frac{1}{n_p!}
=\prod _{p=1}^{dim X}\exp \left(-\frac{t^p \theta_p}{p}\right).
\end{eqnarray*}
\end{proof}

\bigskip
\subsection{Computation of the class S}
We now aim at finding expressions for the class $S$ in terms of the
geometric Segre forms $s'$
$$s'_m(E,h)=\pi_\star(\alpha^{r-1+m})=\left({r-1+m\atop r-1}\right)
\pi_\star \big((-\Theta^1a^\star)^{m}\Omega^{r-1}\big).$$ 
We will need  as auxiliary tools the forms on $X$ defined by
$$s^b_{c}(E,h)
=\pi_\star\big((-1)^b(\Theta^b a^\star)\alpha^{r-1+c}\big)
=\left({r-1+c\atop r-1}\right)
\pi_\star\big((-1)^b(\Theta^b
a^\star)(-\Theta^1a^\star)^c\Omega^{r-1}\big).$$
Integral computations of $\pi_\star$ similar to the previous one lead to
\begin{eqnarray*}
s^b_{c}(E,h)&=&\frac{(-1)^{b+c}}{(r+c)c!}\sum_{I\in\nb_r^{c+1}\atop
  K\in\nb_r^{b-1}}\sum_{\sigma\in\Sigma_{c+1}}c_{i_1i_{\sigma(1)}}
c_{i_2i_{\sigma(2)}}\cdots c_{i_ci_{\sigma(c)}}
c_{i_{c+1}k_2}c_{k_2k_3}\cdots c_{k_bi_{\sigma(c+1)}}
\end{eqnarray*}
The number of permutations in $\Sigma_{c+1}$ of shape
$1^{n_1}2^{n_2}\cdots (c+1)^{n_{c+1}}$ for which $c+1$ is in a cycle
of length $q$ is $\frac{c!}{q^{n_{q}-1}(n_q-1)!\prod_{p\not
  =q}p^{n_p}n_p!}$. The contribution of such a cycle is 
$(-1)^{q+b-1}\theta_{q+b-1}$
for we plug in $\Theta^b a^\star$ instead of $\Theta^1 a^\star$. 
Forward,
\begin{eqnarray*}
s^b_{c}(E,h)
&=&\frac{1}{(r+c)}\sum_{N\in\nb^{dim X}\atop \sum pn_p=c+1}
\sum_{q/n_q\geq 1} \theta_{q+b-1}\left(\frac{\theta_q}{q}\right)^{n_q-1}
 \frac{1}{(n_q-1)!} 
\prod_{p\not =q}\left(\frac{\theta_p}{p}\right)^{n_p}
 \frac{1}{n_p!}\end{eqnarray*}
Changing the order of the summations over $q$ and $N$ and
then summing over $c$, 
\begin{eqnarray*}
\sum_{c\in\nb} t^{b+c}(r+c)s^b_{c}(E,h)&=&\sum_{q=1}^{dim X}
 t^{q+b-1}\theta_{q+b-1}\prod_{p=1}^{dim X}
\exp \left(\frac{t^p\theta_p}{p}\right)\\
&=&\left(\sum_{m\in\nb}t^ms'_m(E,h)\right)\left(\sum_{q=b}^{dim X}
 t^{q}\theta_{q}\right)
\end{eqnarray*}
 which shows in particular that the forms 
$s^b_{c}(E,h)$ are explicitly computable from the geometric Segre 
forms only and that they are closed. 

Define $\mathcal{H}_a^b$ to be $\displaystyle \sum_{i=1}^a \mathcal{H}_{i}
\Big({a\atop i}\Big) \Big({b\atop i}\Big)^{-1}$.
\begin{theo}\label{theo:HH} The secondary form $S$ is given by
\begin{eqnarray*}
S_{m+1}(E,h)&=&-\sum_{a+b+c=m}\mathcal{H}_{r-1-a-b}^{r-1+c}
c_a(E^\star,h)  s^b_{c}(E,h)
\end{eqnarray*}
\end{theo}

\begin{proof}
For degree reason, $S_m=0$ for $m>dim X +1$.
It follows from the relation
 $$c_p(E\ot L,h_E\ot h_L)=\sum_{i+j=p}\left({r-j\atop p-j}\right) 
c_1(L,h_L)^i c_j(E,h_E)$$ for every holomorphic hermitian vector bundle
$(E,h_E)$ and every holomorphic hermitian line bundle $(L,h_L)$ and
from the transgression techniques of~\cite{gs} (Theorem 1.2.2) that
$\widetilde{c_{r}}(\Sigma(1),h)$ is equal to
$\sum_{i+j=r}\alpha ^{i}\widetilde{c_{j}}(\Sigma,h)=
\sum_{i+j=r}\alpha ^{i}\widetilde{c_{j}}$.
Recall that for $j>r$,  $\widetilde{c_{j}}=0$.
\begin{eqnarray*}
S_{m+1}(E,h)&=&
\pi_\star\big(\alpha^{m}\widetilde{c_r}\left(\Sigma(1),h\right)\big) 
 =\sum_{i+j=r+m}\pi_\star\big(\alpha ^{i}\widetilde{c_{j}}\big)\\
&=&\sum_{i+c+j=r+m}\left({i+c\atop c}\right)
\pi_\star\left((\widetilde{c_{j}})_{r-1-i}
\Omega^i(-\Theta^1a^\star)^c\right)  
\end{eqnarray*}
where we used the relation $\alpha=\Omega-\Theta ^1a^\star$ valid on the
whole fiber $\pi^{-1}(x_0)$.

\bigskip
Terms of full relative degree are simpler to compute than the whole
Bott-Chern forms. The proposition below is a weak form of
theorem~\ref{theo} whose full strength is hence not needed for
arithmetic applications.
\begin{prop}\label{prop:full} On the fiber of $\pi$ over $x_0$,
$$(\widetilde{c_{d+1}})_f\Omega^{r-1-f} = -\mathcal{H}_f
\left({r-1\atop f}\right)^{-1}\left({r-1-d+f\atop f} \right)
\sum_{\alpha+\beta=d-f}\pi^\star c_{\alpha}(E^\star)(-1)^\beta
(\Theta^\beta a^\star) \Omega^{r-1}.$$
\end{prop}
\begin{proof} 
Starting from 
\begin{eqnarray*}
\Phi _{d+1}(u)&=&\sum_{2\leq i_1<i_2<\cdots <i_d\leq r}
\sum_{\sigma\in\Sigma_d}\ep (\sigma)\bigwedge _{m=1}^d 
\left(c_{i_m i_{\sigma (m)}} 
+ (1-u)\frac{i}{2\pi} dz_{i_m}\wedge d\zb_{i_{\sigma (m)}}\right)
\end{eqnarray*}
and 
\begin{eqnarray*}
\Omega^{r-1-f}=(r-1-f)!\sum_{A=(a_1<a_2<\cdots <a_{r-1-f})}
\bigwedge_{j=1}^{r-1-f} \frac{i}{2\pi} dz_{a_j}\wedge d\zb_{a_j}
\end{eqnarray*}
we infer (with $I=J\cup A^c$ and $\sigma=\sigma'\sigma''$)
\begin{eqnarray*}
\lefteqn{(\Phi _{d+1}(u))_f\Omega^{r-1-f}}&&\\
&=&(1-u)^f
(r-1-f)!\sum_{|A|=r-1-f} \bigwedge\frac{i}{2\pi} dz_{a_j}\wedge d\zb_{a_j}\\&&
\sum_{2\leq j_1<j_2<\cdots <j_{d-f}\leq r\atop J\subset A}
\sum_{\sigma'\in\Sigma_{d-f}\atop \sigma''\in\Sigma_f}
\ep(\sigma')\bigwedge_{1\leq m\leq d-f} c_{j_m j_{\sigma' (m)}}
\ep(\sigma'')\bigwedge_{k=1\atop A^c=(b_1<b_2<\cdots <b_f)}^f
\frac{i}{2\pi} dz_{b_k}\wedge d\zb_{b_{\sigma''(k)}}\\
&=&(1-u)^f (r-1-f)!f!\frac{\Omega^{r-1}}{(r-1)!}\sum_{A}
\sum_{2\leq j_1<j_2<\cdots <j_{d-f}\leq r\atop J\subset A}
\sum_{\sigma'\in\Sigma_{d-f}}
\ep(\sigma')\bigwedge_{1\leq m\leq d-f} c_{j_m j_{\sigma' (m)}}\\
&=&(1-u)^f
\frac{(r-1-f)!f!}{(r-1)!}\Omega^{r-1}\sum_{|J|=d-f}\sum_{|A|=r-1-f\atop
  A\supset J}det(c'_{J,J})\\
&=&(1-u)^f \frac{(r-1-f)!f!}{(r-1)!}\Omega^{r-1}
\left( {r-1-d+f\atop f}\right) c'_{d-f}(E^\star).
\end{eqnarray*} 
This leads to
$$(\widetilde{c_{d+1}})_f\Omega^{r-1-f} = -\mathcal{H}_f
\left({r-1\atop f}\right)^{-1}\left( {r-1-d+f\atop f}\right)
c'_{d-f}(E^\star
)\Omega^{r-1}$$
at the point $(x_0,[a^\star_0])$.
To get an expression valid on the whole fiber, we refer to
proposition~\ref{prop:c}.
\end{proof} 

Back to the computation of $S_{m+1}$, with $a+b=(j-1)-(r-1-i)=i+j-r$,
\begin{eqnarray*}
\lefteqn{S_{m+1}(E,h) }\\
&=&\!\!-\!\!\sum_{0\leq i\leq r-1\atop a+b+c=m}\!\!\mathcal{H}_{r-1-i}
\left({r-1 \atop i} \right)^{-1}\!\!\left({r-1-a-b \atop r-1-i} \right)
\left({i+c\atop c}\right)
c_a(E^\star,h)
\pi_\star\left((-1)^b \Theta^b
  a^\star(-\Theta^1a^\star)^c\Omega^{r-1}\right) \\ 
&=&\!\!-\!\!\sum_{a+b+c=m}\!\!\!\!\sum_{i=1}^{r-1-a-b}
\mathcal{H}_{i}\frac{(r-1-a-b)!(r-1-i+c)!}{(r-1)!(r-1-i-a-b)!c!}
c_a(E^\star,h)
\pi_\star\left((-1)^b \Theta^b
  a^\star(-\Theta^1a^\star)^c\Omega^{r-1}\right) \\ 
&=&\!\!-\!\!\sum_{a+b+c=m}\!\!\!\!\sum_{i=1}^{r-1-a-b}\mathcal{H}_{i}
\frac{(r-1-a-b)!(r-1+c-i)!}{(r-1-a-b-i)!(r-1+c)!} c_a(E^\star,h)
  s^b_{c}(E,h) \\
&=&\!\!-\!\!\sum_{a+b+c=m}\!\!\!\!\sum_{i=1}^{r-1-a-b}\mathcal{H}_{i}
\left({r-1-a-b \atop i}\right)\left({ r-1+c\atop i}\right)^{-1}
c_a(E^\star,h)  s^b_{c}(E,h).
\end{eqnarray*} ending the proof of theorem~\ref{theo:HH}.
\end{proof}

\bigskip\bigskip
\section{Some arithmetic applications}
\label{arithmetic}

\subsection{On arithmetic characteristic classes}

We complete the scheme proposed by Elkik (\cite{el}, see
 also~\cite{gs2}) for the construction of the arithmetic
 characteristic classes assuming known the construction of the
 arithmetic first Chern class and the push forward operation in
 arithmetic Chow groups. 

\bigskip Let $K$ be a number field, $\mathcal{O}_K$ its ring of
integers and $S:=spec(\mathcal{O}_K)$.
Consider an arithmetic variety $\mathcal{X}$ on $S$ (i.e. a flat
regular projective scheme $\chi$ over $S$ together with the collection of
schemes $\chi_\cb=\coprod_{\sigma :
  K\hookrightarrow\cb}\chi_\sigma$)  and
 an arithmetic vector bundle $\overline{E}$ of rank~$r$ on it (i.e. an
 locally free sheaf of $\mathcal{O}_\chi$-modules together with
 the collection of corresponding vector bundles on
 $\chi_\cb(\cb)$ endowed with a hermitian metric). 
Write $X$ for the smooth variety $\chi_\cb(\cb)$ and $(E,h)$ 
for the induced hermitian vector bundle on $X$.
 All the objects we will consider on $X$ are
 required to be invariant under the conjugaison given by the real
 structure of $\mathcal{X}$. The notation $A^{p,p}(X)$ is now
 used for the subspace of conjugaison invariant forms. 
 We will denote by $\mathcal{D}(X)$ the space of (conjugaison invariant)
 currents on $X$ and by $\widetilde{\mathcal{D}}(X)$ its quotient by
 $Im d' +Im d''$.

\bigskip
We collect some basic facts we will need on arithmetic intersection theory
 (see~\cite{gsihes} for properties of the arithmetic Chow groups and
 ~\cite{gs} for the construction of the arithmetic Chern classes).
The arithmetic Chow groups are defined to be the quotient by the
 subgroup generated by arithmetic principal cycles and pairs of the
 form $(0,\partial u+\overline{\partial} v)$ of the free Abelian group
 on pairs $(Z,g_Z)$  of an algebraic cycle $Z\subset\chi$ and a Green
 current $g_Z\in\widetilde{\mathcal{D}}(X)$ for $Z(\cb )$ defined by
 the requirement that $\delta_{Z(\cb)}+dd^c g_Z$ be the current
 associated with a smooth form in ${A}(X)$. The notation
 $\delta_{Z(\cb)}$ is used for the current of integration along $Z(\cb)$.
There are natural maps
$$\begin{array}{ccccc}
a&:&A^{p-1,p-1}(X)&\to& \widehat{CH}^p(\mathcal{X})\\
 & &\eta&\mapsto& [(0,\eta)]
\end{array}
\ \ \ \ \ \
\begin{array}{ccccc}
\omega&:&\widehat{CH}^p(\mathcal{X})&\to &A^{p,p}(X)\\
& &[(Z,g_Z)]&\mapsto &\delta_{Z(\cb)}+dd^c g_Z
\end{array}$$
The push-forward map on arithmetic Chow groups is defined
component-wise. Hence, $a\pi_\star=\widehat{\pi}_\star a$.
The product in $\widehat{CH}^\star(\mathcal{X})$ is defined by the
formula~: $[(Y,g_Y)]\cdot [(Z,g_Z)]=[(Y\cap Z,g_Y\star g_Z)]$
where the star product is given by
$$g_Y\star g_Z:=g_Y\wedge \delta_Z + \omega([(Y,g_Y)])\wedge
g_Z\in\widetilde{\mathcal{D}}(X).$$
Hence, for $r\in\widehat{CH}^\star(\mathcal{X})$, $r\cdot
a(\eta)=a(\omega(r)\wedge \eta)$.

As for the arithmetic Chern classes, they have the following
properties~:
\begin{itemize}
\item For every arithmetic vector bundle $\overline{E}$ 
on $\mathcal{X}$,
$\omega \big(\widehat{c_t}(\overline{E})\big)=c_t(E,h).$
\item For every exact sequence $(\overline{\mathcal{S}})=
  (0\to\overline S\to\overline{E}\to \overline{Q}\to 0)$
of arithmetic vector bundles on $\mathcal{X}$,
\begin{eqnarray}
\label{chern}
\widehat{c_t}(\overline{E})-\widehat{c_t}(\overline{S})
\widehat{c_t}(\overline{Q})= -a(t\widetilde{c_t}(\mathcal{S}(\cb),h)) 
\end{eqnarray}
\item  For every arithmetic line bundle $\overline L$ and every
  arithmetic vector bundle $\overline{E}$ of rank $r$ on $\mathcal{X}$,
\begin{eqnarray}
\label{l}
\widehat{c_r}(\overline{E}\otimes
\overline{L})=\sum_{p+q=r}\widehat{c_p}(\overline E)
\widehat{c_1}(\overline L)^q.
\end{eqnarray}
\end{itemize}

\bigskip We now explain the construction of the arithmetic
characteristic classes.
Use first the usual hermitian theory on $X$ to define in $A^{d,d}(X)$
the following characteristic forms
\begin{itemize}
\item the Segre forms $s'_d(E,h):=\pi_\star (\Theta
(\mathcal{O}_E(1),h)^{r-1+d})$,
\item the Chern forms $c_d(E,h):=trace(\Lambda^d \Theta (E,h))$ 
\item  and the forms $\theta_d:=trace(\Theta (E,h)^{\otimes
    d})$.
\end{itemize}

Define the generalized Segre forms $s^b_{c}(E,h)\in A^{b+c,b+c}(X)$ by 
$$\sum_{c=0}^{+\infty} t^{b+c}(r+c)s^b_{c}(E,h)=s'_t(E,h)
\sum_{q=b}^{dim X} t^{q}\theta_{q}.$$ 
Consider the (secondary) characteristic forms $S_{m+1}(E,h)$
 and $R_{m+1}(E,h)$ in $A^{m,m}(X)$ given by the relations
  \begin{eqnarray*}
S_{m+1}(E,h)&=&-\sum_{a+b+c=m}\mathcal{H}_{r-1-a-b}^{r-1+c}
c_a(E^\star,h)  s^b_{c}(E,h)\\
R_t(E,h)&=& s'_t(E,h)S_t(E,h).
\end{eqnarray*} 
where $S_t(E,h)=\sum_{m=0}^{+\infty}t^m S_{m+1}(E,h)$ and 
$R_t(E,h)=\sum_{m=0}^{+\infty}t^m R_{m+1}(E,h)$.
Our previous computations ensure that the class of $S_{m+1}(E,h)$ in
$\widetilde{A^{m,m}}(X)$ is
 \begin{eqnarray}\label{sm}
S_{m+1}(E,h)
=\pi_\star(\widetilde{c}_r(\Sigma(1),h)\Theta(\mathcal{O}_E(1),h)^m). 
\end{eqnarray} 

Define the $m$-th geometric Segre class of $\overline{E}$ to be 
 $$\widehat{s'_m}(\overline{E}):=\widehat{\pi}_\star
 \left( \widehat{c_1}(\mathcal{O}_{\overline{E}}(1)) ^{r-1+m} 
\right)$$
in the arithmetic Chow group 
$\widehat{CH}(\mathcal{X})$ of $\mathcal{X}$.
Its $m$-th arithmetic Segre class $\widehat{s_m}(E,h)$ is defined 
 as the class in $\widehat{CH}(\mathcal{X})$ obtained
 by adding in the definition of the $m$-th geometric Segre class
an extra term computed
precisely thanks to the previous top Bott-Chern form~:
$$\widehat{s_m}(\overline{E}):=\widehat{s'_m}(\overline{E})+a(R_m(E,h)).$$
Comparing with the arithmetic Chern class polynomial of $\overline{E^\star}$ as
defined in ~\cite{gs} we get
\begin{theo} In the arithmetic Chow group of $\mathcal{X}$
$$\widehat{c_t}(\overline{E^\star})\cdot\widehat{s_t}(\overline{E})=1.$$
\end{theo}
\begin{proof}
This is in fact a precise form of formula~(\ref{arithmetique})  in the
arithmetic setting.
For the twisted arithmetic Euler sequence $\overline{\Sigma}(1)$, 
$$0\to \mathcal{O}_{\pb (\overline{E}) }\to
\pi^\star\overline{E^\star}\ot\mathcal{O}_{\overline{E}}(1)\to
T_{\pb (\overline{E})/\mathcal{X}}\to 0$$
formulas~(\ref{l}) and~(\ref{chern}) read
\begin{eqnarray*}
\sum_{p+q=r}\widehat{c_p}(\pi^\star\overline{E^\star}) 
\widehat{c_1}(\mathcal{O}_{\overline{E}}(1)) ^q&=&
\widehat{c_r}(\pi^\star\overline{E^\star}\ot\mathcal{O}_{\overline{E}}(1))\\
&=& \widehat{c_1}(\mathcal{O}_{\pb (\overline{E}) }) 
\widehat{c_{r-1}}(T_{\pb (\overline{E})/\mathcal{X}})
-a(\widetilde{c}_r(\Sigma(1),h)) \\
&=&-a(\widetilde{c}_r(\Sigma(1),h))
\end{eqnarray*}
By push-forward $\widehat{\pi}_\star$ 
to the Chow group of $\mathcal{X}$, using the previously recalled
facts and formula~(\ref{sm}), this leads to 
$$\widehat{c_t}(\overline{E^\star})\cdot \widehat{s'_t}(\overline{E})
=1-a(tS_t(E,h))$$ 
that is
$$\widehat{c_t}(\overline{E^\star})\cdot\widehat{s_t}(\overline{E})
=1-a(tS_t(E,h))+\widehat{c_t}(\overline{E^\star}) \cdot a(tR_t(E,h)).$$
But
$$\widehat{c_t}(\overline{E^\star}) \cdot
a(R_t(E,h))=a(\omega(\widehat{c_t}(\overline{E^\star})) R_t(E,h))=
a(c_t(E^\star,h)R_t(E,h))=a(S_t(E,h))$$
thanks to proposition~\ref{prop:closed}.
\end{proof}
This provides an alternative definition of the arithmetic Chern
classes and hence of all the arithmetic characteristic classes
without using the splitting principle. 
This point of view is closer to that of Fulton~\cite{fu}.

\subsection{On the height of $\pb (\overline{E})$}
For any vector bundle $E$ of rank $r$ on a smooth complex compact
analytic manifold $X$ of dimension $n$, we define the analytic height of 
$\pb(E)\stackrel{\pi}\to X$ by 
$$h_{\mathcal{O}_E (1)}(\pb (E))=\int_{\pb (E)} c_1(\mathcal{O}_E
(1))^{n+r-1}.$$ 
It follows from Fubini theorem that 
$$h_{\mathcal{O}_E (1)}(\pb (E))
=\int_X \pi_\star \left(c_1(\mathcal{O}_E (1))^{n+r-1}\right)=\int_X s'_n(E).$$
According to Hartshorne, the ampleness of $E\to X$ is defined to be the
ampleness of $\mathcal{O}_E (1)\to \pb(E)$ hence implies the positivity
of the analytic height of $\pb (E)$ from its very definition.

\bigskip
In the arithmetic setting, for an arithmetic vector bundle
$\overline{E}$ of rank $r$ over an arithmetic variety
$\mathcal{X}\stackrel{f}\to S$  of relative dimension $n$ over $S$,
we define the arithmetic height of
 $\pb(\overline{E})~\stackrel{\pi}\to~\mathcal{X}$ with respect to 
$\mathcal{O}_{\overline{E}}(1)$ by
$$\widehat{h}_{\mathcal{O}_{\overline{E}} (1)}(\pb (\overline{E}))
=\widehat{deg}\widehat{f}_\star\widehat{\pi}_\star
\left(\widehat{c_1}(\mathcal{O}_{\overline{E}} (1))^{dim \pb(\overline{E})}
\right)$$
where $\widehat{deg} : \widehat{CH}^1(S)\to \widehat{CH}^1(\zbb)\to
\rb$ the last map sending $(0,\lambda)$ to $\lambda/2$ (~see~\cite{bgs}).
From the definition of the arithmetic Segre class, we infer
$$\widehat{h}_{\mathcal{O}_{\overline{E}} (1)}(\pb (\overline{E}))
=\widehat{deg} \widehat{f}_\star \widehat{s'}_{n+1} 
=\widehat{deg} \widehat{f}_\star \widehat{s}_{n+1}-\frac{1}{2}\int_X
R_{n+1}.$$ The complex cohomological term $-\frac{1}{2}\int_X
R_{n+1}$ is computed thanks to theorem~\ref{theo:HH}.
For example, 
  \begin{eqnarray*}
    R_1&=&S_1=-\sum_{i=1}^{r-1}\mathcal{H}_i s'_0(E,h)
        =-\sum_{i=1}^{r-1}\mathcal{H}_i .\\
    R_2&=&s'_1S_1+S_2
=-\left(1+\frac{1}{r}\right)\left(\sum_{i=1}^{r-1}\mathcal{H}_i\right)s'_1(E,h)
 \end{eqnarray*}  \begin{eqnarray*}
\textrm{ In rank }2,&&S_{m+1} (E,h)=\frac{1}{m+1}s'_{m}(E,h) \\
             &&R_{m+1} (E,h)=-\sum_{p+q=m}\frac{1}{q+1}s'_p(E,h)s'_{q}(E,h)\\
\textrm{ In rank }3,&& S_{m+1}
        (E,h)=-\frac{1}{m+1}c_1(E^\star,h)s'_{m-1}(E,h)
        -\frac{3m+5}{(m+1)(m+2)}s'_{m}(E,h)\\
&&=\frac{1}{m+1}\left(s'_1(E,h)s'_{m-1}(E,h)-s'_{m}(E,h)\right)
-\frac{2m+3}{(m+1)(m+2)}s'_{m}(E,h).
  \end{eqnarray*}

According to Zhang, a line bundle $\overline{L}\to \mathcal{X}$ is
said to be ample if $L\to\chi$ is ample,  $(L,h)\to X$ is
semi-positive and for every large enough $n$ 
there exists a $\zbb$-basis of $\Gamma (\chi, L^{\otimes  n})$ made of
sections of sup norm less than $1$ on $X$.
This implies that the leading coefficient of the Hilbert function of
$\overline{L}$ is positive (\cite{zhang}, lemma 5.3) which in turn
implies the positivity of the arithmetic $dim \mathcal{X}$-fold
intersection of $\widehat{c_1}(\overline{L})$ is positive (\cite{gs2}).
Hence, the ampleness of $\overline{E}$ that we define to be the
ampleness of the associated line bundle
 $\oc_{\overline{E}}(1)\to\pb(\overline{E})$
implies the positivity of the
arithmetic height of $\pb(\overline{E})$.

\bigskip
 Note however that it does not imply the positivity of the secondary 
term $ -\frac{1}{2}\int_X R_{n+1}$ as it can be checked with
Fulton-Lazarsfeld characterization of numerically positive  
polynomials on ample vector bundles~\cite{fl}~:
 the numerically positive polynomials in Chern classes are non-zero
 polynomials having non-negative coefficients in the basis of Schur
 polynomials in Chern polynomials. By Jacobi-Trudi formula, this basis
 is also the basis of Schur polynomials in Segre classes.
For $r=3$,  $n=3$ the third coefficient of $-\frac{1}{2} R_4$ 
in the degree $3$ part of the basis of Schur polynomials in Segre classes
$(s'_3,s'_1s'_2-s'_3,{s'_3}^3-2s'_1s'_2+s'_3)$
is $-\frac{1}{6}$.

\end{document}